\documentclass[numbers,webpdf,imanum]{ima-authoring-template}
\usepackage[utf8]{inputenc}

\usepackage{amsfonts}
\usepackage{amsmath}
\usepackage{amssymb}
\makeatletter\let\imasavedlsss\l@subsubsection\let\l@subsubsection\@undefined\makeatother
\usepackage{subcaption}
\makeatletter\let\l@subsubsection\imasavedlsss\makeatother
\usepackage{graphicx}

\usepackage{amsthm}
\theoremstyle{thmstylethree} 
\newtheorem{remark}{Remark}[section]
\newtheorem{assumption}{Assumption}[section]

\theoremstyle{thmstyletwo} 
\newtheorem{theorem}{Theorem}[section]
\theoremstyle{thmstyletwo} 
\newtheorem{proposition}{Proposition}[section]
\newtheorem{corollary}{Corollary}[section]

\usepackage{ifpdf}
\usepackage{hyperref}
\newcommand{\email}[1]{\href{mailto:#1}{\texttt{#1}}}

\usepackage{cleveref}
\crefname{remark}{remark}{remarks}
\Crefname{remark}{Remark}{Remarks}

\crefname{assumption}{assumption}{assumptions}
\Crefname{assumption}{Assumption}{Assumptions}

\crefname{theorem}{theorem}{theorems}
\Crefname{theorem}{Theorem}{Theorems}

\crefname{proposition}{proposition}{propositions}
\Crefname{proposition}{Proposition}{Propositions}

\crefname{corollary}{corollary}{corollaries}
\Crefname{corollary}{Corollary}{Corollaries}

\usepackage{algorithm}
\usepackage{algcompatible}
\providecommand{\algorithmicand}{\textbf{and}}\providecommand{\AND}{\algorithmicand{} }
\providecommand{\algorithmicor}{\textbf{or}}\providecommand{\OR}{\algorithmicor{} }

\newcommand{\comm}[1]{\hfill $\triangleright$ \texttt{#1}}
\usepackage{etoolbox}
\makeatletter
\AtBeginEnvironment{algorithmic}{%
  \fontsize{\f@size}{9pt}\selectfont
  \setlength{\abovedisplayskip}{9pt}\setlength{\belowdisplayskip}{8pt}%
  \setlength{\jot}{1pt}}
\makeatother
\makeatletter\providecommand{\theHALG@line}{}\renewcommand{\theHALG@line}{\thealgorithm.\arabic{ALG@line}}\makeatother

\def\blam{\boldsymbol{\lambda}}
\def\bmu{\boldsymbol{\mu}}

\newcommand{\real}{\mathbb{R}}
\newcommand{\feas}{C}

\usepackage{xcolor}

\usepackage{todonotes}

\makeatletter
\def\amsclass#1{\g@addto@macro\@abstract{%
    \par%
    \vbox{\fontsize{9bp}{11}\selectfont\textit{AMS subject classifications:}\ #1}%
}}%
\makeatother

\newif\ifarxiv
\arxivtrue
\ifarxiv
\makeatletter
\patchcmd{\@maketitle}
  {\@journaltitle\rm\ (\@copyrightyear) \textbf{\@vol}, \thepage--\thelastpage}
  {Preprint submitted to \@journaltitle}
  {}{\PackageWarning{arxivswitch}{journal header patch failed}}
\patchcmd{\@maketitle}
  {\href{https://doi.org/\@DOI}{https://doi.org/\@DOI}}
  {}
  {}{\PackageWarning{arxivswitch}{DOI patch failed}}
\patchcmd{\@maketitle}
  {[\@history]}
  {}
  {}{\PackageWarning{arxivswitch}{history patch failed}}
\def\ps@opening{%
  \let\@oddfoot\@empty\let\@evenfoot\@empty
  \let\@evenhead\relax\let\@oddhead\relax}
\makeatother
\fi

\ifpdf
\hypersetup{
  pdftitle={RALES: A Risk-Averse Augmented Lagrangian Equality-constrained SQP Method},
  pdfauthor={H.~Antil, D.~P.~Kouri, D.~Ridzal, and T.~M.~Surowiec}
}
\fi

\begin{document}

\DOI{DOI HERE}
\copyrightyear{2026}
\vol{00}
\pubyear{2026}
\access{Advance Access Publication Date: Day Month Year}
\appnotes{Paper}
\copyrightstatement{Published by Oxford University Press on behalf of the Institute of Mathematics and its Applications. All rights reserved.}
\firstpage{1}

\title[RALES: A Risk-Averse Augmented Lagrangian Equality-Constrained SQP Method]{RALES: A Risk-Averse Augmented Lagrangian Equality-Constrained SQP Method}

\author{
  Harbir~Antil
  \address{Department of Mathematical Sciences and Center for Mathematics and Artificial Intelligence (CMAI), George Mason University, Fairfax, VA 22030, USA.
  (\email{hantil@gmu.edu}, \url{http://math.gmu.edu/\string~hantil/}).}}
\author{
  Drew~P.~Kouri
  \address{Optimization and Uncertainty Quantification, Sandia National Laboratories, Albuquerque, NM 87185, USA. 
  (\email{dpkouri@sandia.gov}, \url{https://www.cs.sandia.gov/cr-dpkouri}).}}
\author{
  Denis~Ridzal
  \address{Optimization and Uncertainty Quantification, Sandia National Laboratories, Albuquerque, NM 87185, USA. 
  (\email{dridzal@sandia.gov}, \url{https://www.sandia.gov/\string~dridzal/}).}} 
\author{
  Thomas~M.~Surowiec
  \address{Department of Scientific Computing and Numerical Analysis, Simula Research Laboratory, Oslo, Norway. (\email{thomasms@simula.no}, \url{https://thomas-surowiec.github.io/}).}}

\authormark{H. Antil, D. P. Kouri, D. Ridzal and T. M. Surowiec}


\received{Date}{0}{Year}
\revised{Date}{0}{Year}
\accepted{Date}{0}{Year}

\abstract{
This paper develops RALES, an augmented Lagrangian method for risk-averse
optimization in function spaces with nonsmooth risk measures, nonlinear
equality constraints, and convex geometric constraints.
The method combines safeguarded risk and constraint multiplier updates
with efficient sequential quadratic programming based on inexact
subproblem solves.
The analysis establishes asymptotic stationarity for nonconvex problems
and gives conditions for feasibility of weak primal accumulation points.
Under additional assumptions on subproblem accuracy, duality, and boundedness,
it provides global primal optimality guarantees and weak convergence of the
full stored and candidate multiplier sequences to a common dual maximizer.
Dual convergence is strong when all penalty parameters remain bounded.
Numerical experiments for a finitely sampled semilinear PDE control problem,
with the Conditional Value-at-Risk and the Mean-Plus-Semideviation risk measures,
illustrate the method’s performance as the spatial mesh is refined and as
the number of samples increases.
}


\amsclass{
49M37, 
90C15, 
90C30, 
90C48, 
49K20, 
49J20, 
65K05  
}

\maketitle

\section{Introduction}
In this paper, we provide a provably convergent algorithm for the 
following class of constrained nonsmooth infinite-dimensional optimization problems:
\begin{equation}\label{eq:opt-prob}
  \min_{x\in X}\; \mathcal{R}(f(x)) \quad\text{subject to}\quad
     g(x) = 0, \quad Tx\in C := \bigcap_{i=1}^m C_i.
\end{equation}
Here, $X$ and $Y$ are real Banach spaces, $Z$ is a real Hilbert space
identified with its topological dual $Z^*$, and
$(\Omega,\mathcal{F},\mathbb{P})$ is a probability space. The mapping $f:X\to
L^2(\Omega,\mathcal{F},\mathbb{P})$ represents an objective function subject to
parametric uncertainty,
$\mathcal{R}:L^2(\Omega,\mathcal{F},\mathbb{P})\to\real$ is a functional
quantifying risk preferences, $g:X\to Y$ is a differentiable, potentially
nonlinear mapping, and $T : X \rightarrow Z$ is a bounded linear operator.
Moreover, the sets $C_1, \dots, C_m$ are nonempty, closed, and convex subsets
of $Z$. 

This formulation is partially motivated by PDE-constrained optimization
with uncertain inputs, risk-averse objectives, and state or control constraints
\cite{Gahururu2022,Geiersbach2022,Geiersbach2021,Henrion2025}.
The nonsmooth risk functional and the geometric constraints require
treatment alongside the nonlinear equations that define the physical
model. In particular, eliminating the equality constraints can require
accurate PDE solves at every objective evaluation.

We introduce RALES, a \textbf{R}isk-Averse \textbf{A}ugmented
\textbf{L}agrangian \textbf{E}quality-Constrained \textbf{S}QP method,
which combines the primal-dual risk minimization
algorithm of \cite{kouri2021primal} with the treatment of explicit
nonlinear equality constraints in ALESQP \cite{Antil2023}. The former
handles nonsmooth risk through a sequence of smooth optimization
problems; the latter uses inexact SQP solves for an augmented Lagrangian
with general geometric constraints. RALES treats the risk identifier
and the geometric-constraint multipliers within a common update scheme,
while controlling equality feasibility and stationarity through the
inner stopping tests. This permits simultaneous updates of states and
controls without requiring exact equality feasibility at each outer
iteration.

The RALES algorithm is a general augmented Lagrangian method that
solves an equality-constrained subproblem at each iteration.  Our
convergence analysis does not depend on how this subproblem is
solved.  As the name suggests, we recommend employing SQP---particularly,
the composite-step approach in \cite{MHeinkenschloss_DRidzal_2014a}.  
In the context of PDE-constrained problems, one could solve the subproblem
using a reduced-space optimization algorithm that eliminates the PDE solution
using the control-to-state map, which may require computationally
expensive nonlinear state solves.
In contrast, SQP applied to the full-space formulation retains both
state and control variables and iterates towards PDE feasibility and optimality
within the same iteration.  Full-space SQP steps therefore make
progress without solving the nonlinear PDE to high accuracy at every
iteration~\cite{Antil2023,MHeinkenschloss_DRidzal_2014a}. This
flexibility, which can accelerate equality-constrained subproblem solves,
comes at the cost of the memory burden associated with larger
coupled linear systems.
For the latter, our implementation uses iterative methods
and preconditioners that exploit the state/control
structure~\cite{rees2010optimal}.

The contributions of this paper are as follows:
\begin{enumerate}
\item We derive a continuously differentiable augmented Lagrangian
for the nonsmooth risk objective and the constraints $Tx\in C_i$.
The risk and constraint multiplier candidates have explicit projection
formulas. These formulas lead to an algorithm with separate penalty
parameters, safeguarded multiplier updates, and equality-constrained
SQP subproblems that can be solved inexactly using matrix-free
methods~\cite{MHeinkenschloss_DRidzal_2014a}.

\item For the nonconvex setting, we analyze the first-order stopping
tests used by the algorithm. For bounded primal sequences, we obtain
vanishing equality residuals and an asymptotic variational inequality
(\Cref{thm:asymstat}). We also give conditions for geometric feasibility
of weak primal accumulation points (\Cref{thm:weak.acc.feas}) and prove
strong convergence of the stored and candidate risk and
geometric-constraint multipliers when all penalties remain bounded
(\Cref{th:conv-dual-finite}).

\item In a separate analysis with exact equality feasibility and
approximate global minimization of the subproblems, we show that weak primal
accumulation points are globally optimal up to the limiting error in the subproblem 
objective value (\Cref{th:conv-primal}). Under the additional minimax
identity, boundedness of the multiplier candidates, and summability of
the inner objective errors weighted by the largest penalty parameter,
the full stored and candidate multiplier sequences
converge weakly to the same dual maximizer, without assuming uniqueness
(\Cref{th:conv-dual-ideal}). This convergence is strong when the
penalties remain bounded.

\item We study the method numerically for a semilinear PDE control
problem with Gaussian forcing, state and control bounds, and two
nonsmooth risk measures. The experiments assess stationarity,
feasibility, and risk identification while varying the spatial mesh
and the number of samples. They show modest growth in the outer
iteration counts over the tested discretizations. The state
constraints in these computations are imposed at finitely many
sampled realizations.
\end{enumerate}

Related work includes stochastic SQP methods for nonlinear programs
with sampled or noisy objective and derivative information
\cite{Berahas2024Stochastic,Oztoprak2023Constrained}. Augmented
Lagrangian methods have also been developed for problems with general
lower-level constraints
\cite{Andreani2008Augmented,EGBirgin_LFBueno_JMMartinez_2016a},
for constrained optimization in Banach spaces
\cite{EBorgens_CKanzow_DSteck_2019a}, and for problems with structured
nonconvex constraints \cite{Jia2022Augmented}. Our main contribution
is the joint treatment of a nonsmooth risk functional and explicit nonlinear
equality constraints, with convergence results stated in function spaces.

Section~\ref{s:not} gives the notation and assumptions.
Sections~\ref{s:aug-lag} and~\ref{s:algo} present the augmented
Lagrangian and the algorithm. Section~\ref{sec:conv} separates the
analysis based on first-order inner accuracy from the stronger results
requiring global subproblem accuracy. Section~\ref{sec:apps} presents
the numerical experiments.

\section{Notation and Assumptions}
\label{s:not}
Throughout, $(X,\|\cdot\|_X)$ and $(Y,\|\cdot\|_Y)$ are real Banach spaces
and $(Z,\|\cdot\|_Z)$ is a real Hilbert space.  We employ standard notation
from convex and functional analysis. For readability, we recall several notions
here. We denote the topological dual space of $X$ by $X^*$ and the associated
duality pairing by $\langle\cdot,\cdot\rangle_{X^*,X}$. We denote by $B_\rho^X$
for $\rho>0$ the closed norm ball on $X$ with radius $\rho$ centered at $0$.
We further denote the space of bounded linear operators that map $X$ into $Y$
by $\mathcal{L}(X,Y)$. For the Hilbert space $Z$, we denote by
$(\cdot,\cdot)_Z$ the inner product on $Z$ and we assume that $\|\cdot\|_Z$ is
the usual Hilbert norm on $Z$. Moreover, for a closed, convex subset $C$ of the
Hilbert space $Z$, we denote the projection of a point $z\in Z$ onto $C$ by
$\mathbf{P}_C(z)$ and the distance from $z$ to $C$ by $d_C(z)$.  That is,
$\mathbf{P}_C(z)$ and $d_C(z)$ satisfy
\[
  d_C(z) := \min_{z'\in C}\;\|z-z'\|_Z = \|z - \mathbf{P}_C(z)\|_Z.
\]
Moreover, we denote the normal cone to $C$ at the point $z\in C$ by
\[
  N_C(z) := \{\lambda\in Z\,\vert\,(\lambda,c-z)_{Z}\le 0
    \;\;\forall\,c\in C\},
\]
with $N_C(z)=\emptyset$ if $z\not\in C$. Finally, we denote convergence with
respect to the weak topology by $\rightharpoonup$, the weak$^*$ topology by
$\rightharpoonup^*$, and the norm topology by $\to$.

Throughout, $(\Omega,\mathcal{F},\mathbb{P})$ is a probability space, where
$\Omega$ is the set of outcomes, $\mathcal{F}\subseteq 2^\Omega$ is a
$\sigma$-algebra of events and $\mathbb{P}:\mathcal{F}\to[0,1]$ is a
probability measure, and $L^2(\Omega,\mathcal{F},\mathbb{P})$ denotes the
(Lebesgue) space of square-integrable random variables defined on
$(\Omega,\mathcal{F},\mathbb{P})$. Moreover,
$f:X\to L^2(\Omega,\mathcal{F},\mathbb{P})$, $g:X\to Y$,
$T\in\mathcal{L}(X,Z)$,
$\mathcal{R}:L^2(\Omega,\mathcal{F},\mathbb{P})\to\real$, and $C\subset Z$
is nonempty, closed and convex.  For the subsequent convergence analysis,
several regularity assumptions are required. We will explicitly state when each
assumption is used. Assumptions (A0) (feasibility), (A1) (differentiability)
and (A2) (risk measure) will be required throughout, whereas (A3), (A4) and
(A5) will only be required to prove convergence.
\begin{assumption}[Regularity of Problem Data]\label{as:diff}
\leavevmode
\begin{enumerate}
  \item[(A0)] There exists $\bar{x}\in X$ such that $g(\bar{x})=0$ and $T\bar{x}\in C$.
  \item[(A1)] The functions $f:X\to L^2(\Omega,\mathcal{F},\mathbb{P})$ and
              $g:X\to Y$ are continuously Fr\'{e}chet differentiable.
  \item[(A2)] The function $\mathcal{R}:L^2(\Omega,\mathcal{F},\mathbb{P}) \to \mathbb R$ is
	      closed, convex, positively homogeneous and monotonic with respect
	      to the pointwise ordering on
              $L^2(\Omega,\mathcal{F},\mathbb{P})$, i.e., for
              $v,\;v'\in L^2(\Omega,\mathcal{F},\mathbb{P})$
              \[
                v \le v' \;\;\text{a.s.} \quad\implies\quad
                \mathcal{R}(v) \le \mathcal{R}(v').
              \]
  \item[(A3)] The composite function $\mathcal{R} \circ f$ 
              is weakly lower semicontinuous.
  \item[(A4)] The adjoint operator $T^*\in\mathcal{L}(Z,X^*)$ is injective.
  \item[(A5)] There exist $\feas_i\subset Z$ for $i=1,\ldots,m$ that are nonempty,
              closed and convex for which $\feas=\feas_1\cap\cdots\cap\feas_m\neq\emptyset$
              and $\{\feas_1,\ldots,\feas_m\}$ is boundedly regular in the sense that
              \[
                \max_{i=1,\ldots,m} d_{\feas_i}(Tx_k)\to 0 \qquad\implies\qquad
                   d_{\feas}(Tx_k) \to 0
              \]
              as $k\to\infty$ for every bounded sequence $\{x_k\}\subset X$.
\end{enumerate}
\end{assumption}

We say that $\bar{x}\in X$ is a (first-order) stationary point of
\eqref{eq:opt-prob} if there exists $\bar{\zeta}\in Y^*$ such that
\begin{equation}\label{eq:fostat}
  -(f'(\bar{x})^*\lambda_0 + g'(\bar{x})^*\bar{\zeta}) \in T^*N_\feas(T\bar{x}),\quad \lambda_0 \in \partial \mathcal{R}(f(\bar{x})),
  \quad\text{and}\quad g(\bar{x})=0.
\end{equation}
Note that this presumes $T\bar{x}\in\feas$ since the normal cone is empty
otherwise.  Of course for \eqref{eq:fostat} to be an optimality condition, we
require additional regularity in the form of a constraint qualification
\cite{ioffe1979theory}.

Before proceeding with the derivation and analysis of the risk-averse augmented
Lagrangian functional, several remarks are in order.
\begin{remark}[Banach Space Valued Constraints]
  As in \cite{EBorgens_CKanzow_DSteck_2019a}, we could consider the case
  where $Z$ is a real Banach space that is densely embedded in a real Hilbert
  space.  However, this would complicate the presentation with little added
  benefit.
\end{remark}

The assumptions on the risk measure in (A2) include three of the four axioms of
the well-known coherent risk measures \cite{Artzner1999}. The axioms here are
known to imply a number of useful properties, which we gather in the following
remark.
\begin{remark}[Assumption~(A2)]\label{rem:A2}
  According to Assumption~(A2), $\mathcal R$ is convex, monotonic and
  finite-valued on $L^2(\Omega,\mathcal{F},\mathbb{P})$. Therefore,
  $\mathcal R$ is continuous and subdifferentiable on
  $L^2(\Omega,\mathcal{F},\mathbb{P})$
  \cite[Proposition~6.6]{AShapiro_DDentcheva_ARuszczynski_2014a}.
  Moreover, \cite[Theorem~6.5]{AShapiro_DDentcheva_ARuszczynski_2014a} ensures
  that there exists a nonempty, closed, and convex set
  $\mathfrak{A}\subseteq\{\theta\in L^2(\Omega,\mathcal{F},\mathbb{P})\,\vert\, \theta \ge 0\;\;\text{a.s.}\}$
  for which
  \[
    \mathcal{R}(v) = \sup_{\theta\in\mathfrak{A}} \; \mathbb{E}[\theta v] \quad\forall\, v\in L^2(\Omega,\mathcal{F},\mathbb{P}).
  \]
  In particular, $\mathfrak{A} = \partial\mathcal{R}(0)$ and therefore,
  $\theta\in\mathfrak{A}$ is bounded by the Lipschitz constant of
  $\mathcal{R}$ near 0.
\end{remark}

Another finer point we highlight here is in regards to (A4).
\begin{remark}[Assumption (A4) as a Constraint Qualification]
  Condition (A4) was used in \cite{Antil2023} to prove asymptotic
  feasibility of iterates of the ALESQP algorithm. As stated in
  \cite[Remark~3.3]{Antil2023}, the regularity assumption (A4) on $T$ is
  \textit{strictly weaker} than Robinson's constraint qualification
  (generalized Mangasarian--Fromovitz) in infinite dimensions. Indeed, recall
  that the operator $T$ is surjective if and only if $T^*$ is injective and the
  range of $T^*$ is norm-closed \cite[Th.~4.15]{rudin1991functional}. Hence,
  (A4) is satisfied if $T$ is surjective.

  Numerous examples of injective bounded linear operators with non-closed
  ranges exist. For example, if $T^*$ is the embedding of the Sobolev space
  $H^1(0,1)$ into $L^2(0,1)$, then $T^*$ is injective, but its range is not
  norm-closed in $L^2(0,1)$.  If it were, then there would exist a constant
  $\alpha > 0$ such that $\|u\|_{L^2}\ge\alpha \|u\|_{H^1}$ for all
  $u\in H^1(0,1)$ \cite[Fact~2.26]{bauschke2017convex}, which is clearly false
  as this would imply that $\|\cdot\|_{L^2}$ and $\|\cdot\|_{H^1}$ are
  equivalent on $H^1(0,1)$.
\end{remark}

\section{A Risk-Averse Augmented Lagrangian Function}
\label{s:aug-lag}

As indicated by its name, RALES builds on the philosophy of augmented
Lagrangian methods by introducing a risk-averse augmented Lagrangian functional
that handles the nonsmooth risk measure $\mathcal R$ alongside the general
constraints $Tx \in C$. We derive this functional here.

We start by defining the indicator functions $I_{\feas_i}$ for each of the
constraint sets $C_i$, i.e., $I_{\feas_i}(y) = 0$ if $y\in{\feas_i}$ and
$I_{\feas_i}(y)=\infty$ if $y\not\in{\feas_i}$. Consequently, we can rewrite
\eqref{eq:opt-prob} as the following nonsmooth optimization problem with
explicit equality constraints:
\begin{equation}\label{eq:opt-prob-2}
  \min_{x\in X}\; \mathcal{R}(f(x)) + \sum_{i=1}^m I_{\feas_i}(Tx) \quad\text{subject to}\quad
     g(x) = 0.
\end{equation}
Since each $\feas_i$ is a nonempty, closed, and convex set, the Fenchel-Moreau
theorem guarantees that $I_{\feas_i} = I^{**}_{\feas_i}$, where
\[
I^{**}_{\feas_i}(y) = \sup_{\lambda_i \in Z}\left\{ (\lambda_i, y)_{Z} - I^*_{\feas_i}(\lambda_i) \right\}.
\]
Coupling this with the observation in \Cref{rem:A2}, we can rewrite
\eqref{eq:opt-prob-2} (and therefore \eqref{eq:opt-prob}) as the min-max
problem:
\begin{equation}\label{eq:opt-prob-1}
  \min_{x\in X}\sup_{\substack{\lambda_0\in\mathfrak{A} \\ \lambda_i\in Z}}\;
    \left\{ \mathbb{E}[\lambda_0 f(x)] + \sum_{i=1}^m (\lambda_i,Tx)_Z - I^*_{\feas_i}(\lambda_i) \right\}
    \quad\text{subject to}\quad
    g(x) = 0.
\end{equation}
Due to its similarity to the smooth setting, we will refer to the inner
objective functional as the Lagrangian.

To simplify the presentation, we denote the (Hilbert) space of dual variables
by $\mathfrak{Z}$ and the set of admissible dual variables by $\mathfrak{L}$,
i.e.,
\[
\mathfrak{Z} := L^2(\Omega,\mathcal{F},\mathbb P) \times \prod_{i=1}^m Z\qquad \text{ and }\qquad
\mathfrak{L} := \mathfrak{A} \times \prod_{i=1}^m Z.
\]
Moreover, we define the vector of all dual variables $\blam$ and the subvector
$\blam_{\bullet}$ for the constraints by 
\[
  \blam := (\lambda_0, \lambda_1, \dots, \lambda_m) \in \mathfrak{Z} \qquad \text{ and } \qquad
  \blam_{\bullet} := (\lambda_1, \dots, \lambda_m) \in \prod_{i=1}^m Z.
\]
We split the Lagrangian, $\ell:X\times\mathfrak{Z}\to\real\cup\{-\infty\}$, for
\eqref{eq:opt-prob-1} into two parts
\[
\ell(x,\boldsymbol{\lambda}) := \ell^1(x,\lambda_0) + \ell^2(x,\boldsymbol{\lambda}_{\bullet}),
\]
where
$\ell^1: X\times \mathfrak{A} \to \mathbb R$ and
$\ell^2: X \times \prod_{i=1}^m Z \to \mathbb R\cup\{-\infty\}$:
\[
\ell^1(x,\lambda_0) := \mathbb E[\lambda_0 f(x)] \qquad \text{and} \qquad
\ell^2(x,\blam_{\bullet}) := \sum_{i=1}^m \left( (\lambda_i, Tx)_{Z} - I^*_{C_i}(\lambda_i) \right).
\]

Given a candidate vector of dual variables $\boldsymbol{\lambda}$ and positive
numbers $r_i>0$ for $i=0,1,\ldots,m$, we obtain a relaxed, augmented Lagrangian
subproblem by subtracting quadratic penalties from the min-max objective
function in \eqref{eq:opt-prob-1} (cf.\ \cite{RTRockafellar_1976b}), yielding
\begin{equation}\label{eq:opt-prob-1r}
  \begin{aligned}
  &\min_{x\in X}\sup_{\boldsymbol{\mu}\in\mathfrak{L}} 
     \left\{ \mathbb{E}[\mu_0 f(x)] - \tfrac{1}{2r_0}\mathbb{E}[(\mu_0-\lambda_0)^2]
    + \sum_{i=1}^m \left[(\mu_i,Tx)_Z - I^*_{\feas_i}(\mu_i) - \tfrac{1}{2r_i}\|\mu_i-\lambda_i\|_Z^2\right] \right\} \\
  &\text{subject to}\quad g(x) = 0.
  \end{aligned}
\end{equation}
The objective function in \eqref{eq:opt-prob-1r} is the risk-averse augmented
Lagrangian functional
\begin{equation}\label{eq:aug-lag}
  L(x,\blam,\mathbf{r})
   := \Psi_0(x,\lambda_0,r_0) + \sum_{i=1}^m\Psi_i(x,\lambda_i,r_i),
\end{equation}
where $\mathbf{r}=(r_0,r_1,\ldots,r_m)\in\real_{>0}^{m+1}$,
$\Psi_0(x,\lambda_0,r_0):=\overline{\Psi}_0(f(x),\lambda_0,r_0)$ with
$\overline{\Psi}_0:L^2(\Omega,\mathcal{F},\mathbb{P})\times L^2(\Omega,\mathcal{F},\mathbb{P})\times(0,\infty)\to\real$
is
\begin{equation}\label{eq:pen0}
  \overline{\Psi}_0(v,\lambda,r) := \sup_{\mu\in \mathfrak{A}}\{\mathbb{E}[\mu v] - \tfrac{1}{2r}\mathbb{E}[(\mu-\lambda)^2]\}
\end{equation}
with $\lambda_0\in\mathfrak{A}$ and $r_0>0$, and
$\Psi_i(x,\lambda_i,r_i):=\overline{\Psi}_i(Tx,\lambda_i,r_i)$ with
$\overline{\Psi}_i:Z\times Z\times(0,\infty)\to\real$ for $i=1,\ldots,m$ is
\begin{equation}\label{eq:pen}
  \overline{\Psi}_i(z,\lambda,r) := \sup_{\mu\in Z}\left\{(\mu,z)_Z - I_{\feas_i}^*(\mu) - \tfrac{1}{2r}\|\mu-\lambda\|_Z^2\right\},
\end{equation}
with $\lambda\in Z$ and $r>0$ for $i=1,\ldots,m$.

To summarize the explicit forms and differentiability of these components for
subsequent analysis, we gather the following properties into a conveniently
citable result.

\begin{proposition}\label{prop:penalty_props}
For any $r>0$, the penalty functionals $\overline{\Psi}_0(\cdot,\lambda,r)$ and
$\overline{\Psi}_i(\cdot,\lambda,r)$ are continuous, convex, and Fr\'{e}chet
differentiable. Their explicit forms satisfy
\begin{subequations}
\begin{align}
  \overline{\Psi}_0(v,\lambda,r) &= \frac{r}{2} \mathbb{E}[v^2] + \mathbb{E}[\lambda v] - \frac{1}{2r}\mathbb{E}\left[|({\rm Id} - \mathbf{P}_{\mathfrak{A}})(rv + \lambda)|^2\right] \quad \lambda \in L^2(\Omega,\mathcal{F},\mathbb P), \label{eq:psi0_explicit} \\
  \overline{\Psi}_i(z,\lambda,r) &= \frac{1}{2r}\|r({\rm Id} - \mathbf{P}_{C_i})(r^{-1} \lambda + z)\|^2_{Z} - \frac{1}{2r} \|\lambda\|^2_{Z}\quad \lambda \in Z. \label{eq:psii_explicit}
\end{align}
\end{subequations}
Here, ${\rm Id}$ denotes the identity map on either
$L^2(\Omega,\mathcal{F},\mathbb{P})$ or $Z$ depending on the context.
Furthermore, the suprema in \eqref{eq:pen0} and \eqref{eq:pen} are attained by
the unique maximizers
\begin{subequations}
\begin{align}
  \overline{\Lambda}_0(v,\lambda,r) &:= \mathbf{P}_{\mathfrak{A}}(rv + \lambda), \label{eq:aug-lag-max-0} \\
  \overline{\Lambda}_i(z,\lambda,r) &:= r((r^{-1}\lambda+z) - \mathbf{P}_{\feas_i}(r^{-1}\lambda + z)), \label{eq:aug-lag-max}
\end{align}
\end{subequations}
and the respective gradients of $\overline{\Psi}_0(\cdot,\lambda,r)$ and
$\overline{\Psi}_i(\cdot,\lambda,r)$ are given by 
\begin{subequations}
\begin{align}
  \nabla_v \overline{\Psi}_0(v,\lambda,r) &= \overline{\Lambda}_0(v,\lambda,r), \\
  \nabla_z \overline{\Psi}_i(z,\lambda,r) &= \overline{\Lambda}_i(z,\lambda,r).
\end{align}
\end{subequations}
Finally, the gradients $\nabla_v\overline{\Psi}_0(\cdot,\lambda,r)$ and
$\nabla_z\overline{\Psi}_i(\cdot,\lambda,r)$ for $i=1,\ldots,m$ are Lipschitz continuous.
\end{proposition}

\begin{proof}
By simple algebraic manipulations and the definition of the projection in
$L^2(\Omega,\mathcal{F},\mathbb{P})$ onto the closed convex set	$\mathfrak{A}$,
we see that 
\[
\aligned
\overline\Psi_0(v,\lambda,r) 
&= 
\frac{r}{2} \mathbb E[v^2] + \mathbb E[\lambda v] - \frac{1}{r} \inf_{\mu \in \mathfrak{A}}\left\{ \frac{1}{2} \| \mu - (r v + \lambda) \|^2_{L^2(\Omega,\mathcal{F},\mathbb{P})}\right\}\\
&=
\frac{r}{2} \mathbb E[v^2] + \mathbb E[\lambda v]  - \frac{1}{2r}\mathbb E\left[|({\rm Id} - \mathbf{P}_{\mathfrak{A}})(rv + \lambda)|^2\right].
\endaligned
\]
We denote the unique minimizer for the infimum above by
$\overline\Lambda_0(v,\lambda,r)$ as defined in \eqref{eq:aug-lag-max-0}.
Although $\overline\Psi_0$ contains a projection operator, which is in general
nonsmooth, one can show that $\overline\Psi_0(\cdot,\lambda,r)$ is the infimal
convolution of a proper, lower semicontinuous, convex function and a smooth
strongly convex function and is therefore continuous, convex, and Fr\'echet
differentiable with a Lipschitz continuous gradient
\cite[Sec.~3]{kouri2021primal}.
Similar arguments apply to $\overline\Psi_i(\cdot,\lambda,r)$, $i=1,\ldots,m$.  See
\cite[Prop.~4.3]{Antil2023} for details.
\end{proof}

A consequence of \Cref{prop:penalty_props} is that the augmented
Lagrangian functional $L$ is continuously Fr\'{e}chet differentiable with
respect to $x$, and its derivative, denoted $L_x'$, is given by
\[
  L_x'(x,\blam,\mathbf r) 
   = \mathbb{E}[\overline\Lambda_0(f(x),\lambda_0,r_0)\nabla f(x)] + \sum_{i=1}^m T^*\overline\Lambda_i(Tx,\lambda_i,r_i).
\]
To simplify the presentation of the RALES algorithm, we denote
\[
  \Lambda_0(x,\lambda_0,r_0):=\overline\Lambda_0(f(x),\lambda_0,r_0)
  \quad\text{and}\quad
  \Lambda_i(x,\lambda_i,r_i):=\overline\Lambda_i(Tx,\lambda_i,r_i)
  \quad\text{for $i=1,\ldots,m$.}
\]
Moreover, we denote by $\|\cdot\|_i$ the norm on the space associated with
$\Lambda_i$, i.e., $\|\cdot\|_0$ is the $L^2(\Omega,\mathcal{F},\mathbb{P})$
norm, while $\|\cdot\|_i$ is the $Z$ norm for $i=1,\ldots,m$.

\section{The RALES Algorithm}
\label{s:algo}
We are now in a position to derive and analyze the RALES algorithm. Motivated
by
\cite{Antil2023,EGBirgin_LFBueno_JMMartinez_2016a,ARConn_NIMGould_PhToint_1991a,kouri2021primal},
we formulate \Cref{algo:rales-method} using the augmented Lagrangian
$L$. Due to the large number of parameters and use of inexact error control and
parameter updating steps, we present the method in three parts: a high-level
overview (\Cref{algo:rales-method}), a collection of all the
initialization parameters (\Cref{algo:rales-init}), and detailed
implementation specification  (\Cref{algo:rales-update}). At its core, RALES
alternates between an inexact primal phase and a dual update phase. In the
primal phase, we seek an approximate stationary point of the augmented
Lagrangian $L(\cdot,\blam,\mathbf{r})$ subject to the equality constraints
$g(x) = 0$. In practice, this step is driven by a composite-step SQP
trust-region method that linearizes the equality constraints. Once an
acceptable primal step is found, the dual phase updates the risk identifier
$\lambda_0$, the constraint multipliers $\lambda_i$, and the penalty parameters
$\mathbf{r}$ based on the measured progress toward feasibility.

\begin{algorithm}[!h]
\caption{RALES Algorithm}
\label{algo:rales-method}
\small
\begin{algorithmic}[1]
\REQUIRE{Initialize using \Cref{algo:rales-init}.}
\FOR{$k=1, 2, 3, \ldots$}
\STATE \textbf{Subproblem Solve:} Compute $(x^{(k)}, \zeta^{(k)}) \in X \times Y^*$
       satisfying
       \vspace{-4pt}
       \[
         \|L'_x(x^{(k)},\blam^{(k)},\mathbf{r}^{(k)})+g'(x^{(k)})^*\zeta^{(k)}\|_{X^*}\le\varepsilon^{(k)}
         \qquad\text{and}\qquad \|g(x^{(k)})\|_Y\le \delta^{(k)}.
       \]
       \vspace{-14pt}
\STATE \textbf{Termination Check:} Exit if $x^{(k)}$ satisfies
       \vspace{-2pt}
       \[
         \begin{aligned}
         &\|L'_x(x^{(k)},\blam^{(k)},\mathbf{r}^{(k)})+g'(x^{(k)})^*\zeta^{(k)}\|_{X^*}\le\varepsilon_*,
	 \qquad \|g(x^{(k)})\|_Y\le\delta_*, \\
         &\max_{1\le i\le m}
         \|\Lambda_i(x^{(k)},\lambda_i^{(k)},r_i^{(k)})
                    -\lambda_i^{(k)}\|_i \le r_i^{(k)} \tau_*, \\
         &\text{and}\qquad \|\lambda_0^{(k)}-\Lambda_0(x^{(k)},\lambda_0^{(k)},r_0^{(k)})\|_0
	 \le r_0^{(k)}\sigma_*.
	 \end{aligned}
       \]
       \vspace{-8pt}
\STATE \textbf{Update Dual Variables and Penalty Parameters:} Update $r_i^{(k+1)}$ and
       $\lambda_i^{(k+1)}$ using \Cref{algo:rales-update}.
\ENDFOR
\end{algorithmic}
\end{algorithm}

The multiplier stopping test controls infeasibility and ensures that each 
candidate multiplier is normal to its constraint set at a point within 
$\tau_*\) of \(Tx^{(k)}$.
In particular, for $i=1,\ldots,m$, set
\[
 y_i^{(k)}:=\mathbf P_{C_i}
 \left(Tx^{(k)}+(r_i^{(k)})^{-1}\lambda_i^{(k)}\right),
 \qquad
 \widehat\lambda_i^{(k)}:=\Lambda_i(x^{(k)},\lambda_i^{(k)},r_i^{(k)}).
\]
The projection formula gives
$\widehat\lambda_i^{(k)}\in N_{C_i}(y_i^{(k)})$ and
\[
 d_{C_i}(Tx^{(k)})\le\|Tx^{(k)}-y_i^{(k)}\|_Z
 =\frac{1}{r_i^{(k)}}\|\widehat\lambda_i^{(k)}-\lambda_i^{(k)}\|_i.
\]
If all stopping residuals vanish, then $y_i^{(k)}=Tx^{(k)}$.
The projection condition for $i = 0$, involving the nonsmooth risk measure, gives
$\lambda_0^{(k)}\in\partial\mathcal R(f(x^{(k)}))$.
Together with stationarity and feasibility of the equality constraint, this implies
\eqref{eq:fostat}, since
$\sum_i N_{C_i}(Tx^{(k)})\subseteq N_C(Tx^{(k)})$.
At positive tolerances, the multiplier test provides a point $y_i^{(k)}$
within $\tau_*$ of $Tx^{(k)}$ at which $\widehat\lambda_i^{(k)}$ is normal.

\begin{algorithm}[!ht]
  \caption{RALES Initialization and Parameters}
  \label{algo:rales-init}
  \renewcommand{\algorithmicensure}{\textbf{Compute:}}
  \small
  \begin{algorithmic}[1]
  \REQUIRE \textbf{Initial Estimates:}
           \STATE $\blam_1 \in \mathfrak{L}$ \comm{Initial dual estimate}
           \STATE $\mathbf{r}_1 \in \mathbb R_{> 0}^{m+1}$ \comm{Vector of positive penalty parameters}
  
  \REQUIRE \textbf{Tolerances and Constants:}
	   \STATE $\{\varepsilon^{(k)}\}$, $\{\delta^{(k)}\} \subset \mathbb R_{>0}$ with $\varepsilon^{(k)}\to0$, $\delta^{(k)}\to0$ \comm{Primal and infeasibility tolerances}
           \STATE $K_0 \in \mathbb N$ \comm{Aggregation index for penalty parameters}
           \STATE $\boldsymbol{\eta} \in \mathbb R^{m+1}_{>0}$ with $\eta_i > 1$ \comm{Pre-$K_0$ scaling factors for penalty updates}
           \STATE $\bar{\eta} > 1$ \comm{Post-$K_0$ scaling factors for penalty updates}
	   \STATE $\boldsymbol{\gamma} \in \mathbb R^{m+1}_{>0}$ with $\gamma_i < 1/2$  \comm{Power factors for penalty updates}
	   \STATE $\boldsymbol{\alpha}, \boldsymbol{\beta} \in \mathbb R^{m+1}_{>0}$ \comm{Powers for penalty updates}
	   \STATE $\boldsymbol{\nu}\in\mathbb{R}^{m+1}_{>0}$ \comm{Multiplier safeguard constants}
	   \STATE $\boldsymbol{\tau}^{(0)} \in \mathbb R^{m+1}_{> 0}$ \comm{Dual tolerances}
           \STATE $\boldsymbol{\theta} \in \mathbb R^{m+1}_{> 0}$ with $\theta_i < 1$ \comm{Dual tolerance forcing factors}
           \STATE $\sigma_*, \tau_*, \delta_*, \varepsilon_* > 0$ \comm{Stopping: risk, multiplier, $g(x) = 0$, stationarity}

  \ENSURE \textbf{Initial Forcing Parameters:}
  \STATE Set $\theta_i^{(1)} \gets \min\{1/r_i^{(1)}, \theta_i\}$ and $\tau_i^{(1)} \gets \tau_i^{(0)} (\theta_i^{(1)})^{\alpha_i}$ for $i=0,\dots,m$.
  
  \end{algorithmic}
\end{algorithm}

\begin{algorithm}[!ht]
  \caption{Multiplier and Penalty Update}
  \label{algo:rales-update}
  \small
  \begin{algorithmic}[1]
    \IF{$k = K_0 + 1$}
      \STATE $\eta_i \gets \overline{\eta}$ for $i=0,1,\ldots,m$
    \ENDIF 

    \STATE \texttt{update} $\gets$ \textbf{false}
    \IF{$k > K_0$ \AND $\exists i \in \{0,1,\dots,m\}$ such that $\|\Lambda_i(x^{(k)}, \lambda_i^{(k)}, r_i^{(k)}) - \lambda_i^{(k)}\|_i > r_i^{(k)} \tau_i^{(k)}$}
      \STATE \texttt{update} $\gets$ \textbf{true}
    \ENDIF

    \FOR{$i=0,1, \dots, m$}
      \IF{$\|\Lambda_i(x^{(k)}, \lambda_i^{(k)}, r_i^{(k)}) - \lambda_i^{(k)}\|_i > r_i^{(k)} \tau_i^{(k)}$ \OR \texttt{update} = \textbf{true}}
        \STATE $r_i^{(k+1)} \gets \eta_i r_i^{(k)}$
        \STATE $\theta_i^{(k+1)} \gets \min\{1/r_i^{(k+1)}, \theta_i\}$
        \STATE $\tau_i^{(k+1)} \gets \tau_i^{(0)} (\theta_i^{(k+1)})^{\alpha_i}$
      \ELSE
        \STATE $r_i^{(k+1)} \gets r_i^{(k)}$
        \STATE $\theta_i^{(k+1)} \gets \min\{1/r_i^{(k+1)}, \theta_i\}$
        \STATE $\tau_i^{(k+1)} \gets \tau_i^{(k)} (\theta_i^{(k+1)})^{\beta_i}$
      \ENDIF
      
      \IF{$\|\Lambda_i(x^{(k)}, \lambda_i^{(k)}, r_i^{(k)})\|_i \le \nu_i (r_i^{(k+1)})^{\gamma_i}$}
        \STATE $\lambda_i^{(k+1)} \gets \Lambda_i(x^{(k)}, \lambda_i^{(k)}, r_i^{(k)})$
      \ELSE
        \STATE $\lambda_i^{(k+1)} \gets \lambda_i^{(k)}$
      \ENDIF
    \ENDFOR
  \end{algorithmic}
\end{algorithm}

\begin{remark}[Penalty Parameter and Multiplier Update] \label{rem:penalty-update}
  The penalty parameter update in \Cref{algo:rales-update} is completely
  decoupled for the first $K_0$ iterations. Here, $K_0$ can be taken
  arbitrarily large, but finite, e.g., $K_0=1000$. This allows each $r_i^{(k)}$
  to be calibrated to the scaling associated with the $i$th constraint as well
  as the risk measure ($i=0$). 
  The penalty parameter and multiplier update schemes are adapted from
  \cite[Alg.~2]{ARConn_NIMGould_PhToint_1991a}. In particular, $\nu_i$ and
  $\gamma_i$ are used to ensure that the $i$th multiplier update does not grow
  too rapidly while $\eta_i$ is used to increase the penalty parameter $r_i$.
  After $K_0$ iterations, \Cref{algo:rales-update} switches schemes
  and updates all penalty parameters in unison. This penalty update scheme is a
  safeguard for the case in which \Cref{algo:rales-method} produces an infinite
  sequence of iterations, forcing the sequence to accumulate at a feasible
  point (under certain assumptions), and is typically never active in practice.
\end{remark}

\begin{remark}[Subproblem Tolerance Sequences] \label{rem:tolerance-update}
  In Algorithm~\ref{algo:rales-method}, we update the feasibility tolerances
  $\tau_i^{(k)}$ for $i=0,1,\ldots,m$ in an analogous fashion as
  \cite[Alg.~2]{ARConn_NIMGould_PhToint_1991a} using the parameters
  $\alpha_i$, $\beta_i$, and $\theta_i$. In contrast, the sequences of
  optimality tolerances $\{\varepsilon^{(k)}\}$ and equality constraint
  feasibility tolerances $\{\delta^{(k)}\}$, are only required to be
  nonnegative and converge to zero. A basic choice is
  \begin{align*}
    \delta^{(k+1)} &= \left\{\begin{array}{ll}
      \delta_0\delta^{(k)} & \text{if $\exists\,i$ such that $r_i^{(k+1)}\neq r_i^{(k)}$} \\
      \delta_1\delta^{(k)} & \text{otherwise}
    \end{array}\right.
    \quad\text{and}\quad \\
    \varepsilon^{(k+1)} &= \left\{\begin{array}{ll}
      \varepsilon_0\varepsilon^{(k)} & \text{if $\exists\,i$ such that $r_i^{(k+1)}\neq r_i^{(k)}$} \\
      \varepsilon_1\varepsilon^{(k)} & \text{otherwise}
    \end{array}\right.
  \end{align*}
  for $0<\delta_1 < \delta_0 < 1$ and $0 < \varepsilon_1 < \varepsilon_0 < 1$.
  More complicated updating strategies are possible. For example, adapting the
  tolerance sequence from \cite{ARConn_NIMGould_PhToint_1991a} yields
  \begin{align*}
    \delta^{(k+1)} &= \left\{\begin{array}{ll}
       \delta^{(0)}(\theta^{(k+1)})^{\alpha_\delta} & \text{if $\exists\,i$ such that $r_i^{(k+1)}\neq r_i^{(k)}$} \\
       \delta^{(k)}(\theta^{(k+1)})^{\beta_\delta}  & \text{otherwise}
    \end{array}\right.
    \;\text{and}\quad \\
    \varepsilon^{(k+1)} &= \left\{\begin{array}{ll}
       \varepsilon^{(0)}(\theta^{(k+1)})^{\alpha_\varepsilon} & \text{if $\exists\,i$ such that $r_i^{(k+1)}\neq r_i^{(k)}$} \\
       \varepsilon^{(k)}(\theta^{(k+1)})^{\beta_\varepsilon}  & \text{otherwise}
    \end{array}\right.
  \end{align*}
  for positive constants $\delta^{(0)}$, $\varepsilon^{(0)}$, $\alpha_\delta$,
  $\alpha_\varepsilon$, $\beta_\delta$ and $\beta_\varepsilon$, where, e.g.,
  $\theta^{(k+1)} = (\max_i r_i^{(k+1)})^{-1}$, whereas adapting the sequences
  from \cite{EGBirgin_LFBueno_JMMartinez_2016a} yields
  \begin{align*}
    \delta^{(k+1)}&=\min\{\delta_0\delta^{(k)},\delta_1\|g(x^{(k)})\|_Y\}
    \;\;\text{and}\;\; \\
    \varepsilon^{(k+1)}&=\min\{\varepsilon_0\varepsilon^{(k)},
    \varepsilon_1\|f'(x^{(k)})+g'(x^{(k)})^*\zeta^{(k)}\|_{X^*}\}
  \end{align*}
  for constants $\delta_0,\,\delta_1,\,\varepsilon_0,\,\varepsilon_1\in(0,1)$.
\end{remark}

\begin{remark}[The SQP Subproblem Solve] \label{rem:sqp-subproblem}
  The subproblem solve in Step~2 of Algorithm~\ref{algo:rales-method} abstracts
  away the explicit sequential quadratic programming machinery. In our
  practical implementation, this inexact minimization step is executed
  using the composite-step trust-region SQP algorithm in
  \cite{MHeinkenschloss_DRidzal_2014a}. This approach computes a quasi-normal
  step to improve feasibility of the nonlinear equality constraints $g(x) = 0$
  followed by a tangential step in the null space of the linearized constraint
  Jacobian to minimize a quadratic model of the augmented Lagrangian $L$.  In
  the context of PDE-constrained optimization, 
  the nonlinear state
  equation is satisfied progressively through these steps; a fully
  converged nonlinear state solve is not required at each SQP iteration.
\end{remark}

\section{Convergence Theory}\label{sec:conv}
The convergence theory of \Cref{algo:rales-method} requires several steps for
both the primal variables $x$ and the dual variables $\blam$. Fortunately, the
results in \cite{kouri2021primal} and \cite{Antil2023} provide blueprints for
the individual components, as \Cref{algo:rales-method} combines features of the
PD-Risk algorithm \cite{kouri2021primal} with the ALESQP algorithm
\cite{Antil2023}. 

To better align our theory with practical computational realities, we structure
our analysis in two parts. First, in \Cref{ssec:conv-stat}, we analyze the
asymptotic stationarity and dual convergence of sequences generated by the
algorithm under realistic inexact tolerances. Subsequently, in
\Cref{ssec:conv-dual}, we analyze the global convergence of minimizers in an
idealized setting where the equality constraint $g(x)=0$ can be satisfied
exactly at each iteration. This corresponds to the situation where $\delta^{(k)} =
\delta^* = 0$ throughout all iterations and amounts to an implicit treatment of
the equality constraint, as is often done in PDE-constrained optimization.

\subsection{Asymptotic Stationarity and Dual Convergence}\label{ssec:conv-stat}
In this section, we provide a convergence theory for stationary points with
minimal regularity assumptions. We will see that much of the analysis from
\cite{Antil2023} can be carried over, despite the introduction of the nonsmooth
objective function $\mathcal{R}\circ f$, which was not considered there.

Due to the nature of the adaptive tolerance strategies employed in RALES, there
are several possible ways for the algorithm to terminate. Throughout this
subsection, we will work with sequences generated by \Cref{algo:rales-method}
under the assumption that 
\[
\varepsilon_{*} = \delta_{*} = \tau_{*} = \sigma_{*} = 0.
\] 
Under this standing assumption, we define two types of index sets:
\begin{align*}
\intertext{\noindent \textbf{Penalty Update:} Add $k$ to $\mathcal P_i \subset \mathbb N$ if\vspace{-7pt}}
&\| \Lambda_i(x^{(k)},\lambda^{(k)}_i,r^{(k)}_i) - \lambda^{(k)}_i \|_i > r^{(k)}_i \tau^{(k)}_i, & \text{(Increase penalty parameter } r^{(k)}_i\text{)} \\[-9pt]
\intertext{\noindent\textbf{Dual Update:} Add $k$ to $\mathcal M_i\subset \mathbb N$ if \vspace{-7pt}}
&\|\Lambda_i(x^{(k)}, \lambda_i^{(k)}, r_i^{(k)})\|_i \le \nu_i (r_i^{(k+1)})^{\gamma_i}, & \text{(Update approximate dual variable } \lambda^{(k+1)}_i\text{).} \\[-9pt] 
\end{align*}
In what follows, we denote the complement of a set $A \subset \mathbb N$ by
$A^c := \mathbb N \setminus A$ and the cardinality of $A$ by $|A|$. 

\begin{remark}[Penalty Parameter Asymptotics]
The update strategies in RALES are designed so that:
\[
\exists\,i\in\{0,1,\ldots,m\} \quad\text{such that}\quad r^{(k)}_i \to +\infty \quad\iff\quad \exists\,j \in \{0,1,\dots,m\}\quad\text{such that}\quad |\mathcal P_j| = +\infty.
\]
To see this, note that $r^{(k)}_i \to +\infty$ implies that the update
$r^{(k+1)}_i = \eta_i r^{(k)}_i$ must appear infinitely often, which only happens
if $\|\Lambda_{i}(x^{(k)},\lambda^{(k)}_i,r^{(k)}_i) - \lambda^{(k)}_i \|_i > r^{(k)}_i\tau^{(k)}_i$
or $\tt{update = true}$. If
$\|\Lambda_{i}(x^{(k)},\lambda^{(k)}_i,r^{(k)}_i) - \lambda^{(k)}_i \|_i > r^{(k)}_i\tau^{(k)}_i$
happens infinitely often, then $j = i$, i.e., $|\mathcal P_i|= +\infty$.
Otherwise, we eventually have $k \ge K_0$ and
$\|\Lambda_{j}(x^{(k)},\lambda^{(k)}_j,r^{(k)}_j) - \lambda^{(k)}_j \|_j > r^{(k)}_j\tau^{(k)}_j$ for
some $j \ne i$. The latter condition must occur infinitely often, otherwise
$r^{(k)}_i$ will not be increased. Since there are finitely many constraints $m$,
$|\mathcal P_j| = +\infty$ must hold for some $j \in \{0,1,\dots,m\}$. Conversely,
if $|\mathcal P_j| =+\infty$ for some $j$, then $\tt{update = true}$ holds
infinitely often, in which case we increase all $r^{(k)}_i$ to $+\infty.$
\end{remark}

\begin{theorem}[Strong Dual Convergence]\label{th:conv-dual-finite}
Suppose Assumptions~(A0)-(A2) hold. Let
$\{(x^{(k)},\blam^{(k)},\mathbf r^{(k)})\}$ be an infinite sequence generated
by \Cref{algo:rales-method}, with the stopping conditions ignored and define
the candidate multipliers by
\[
  \widehat\lambda_i^{(k)} :=\Lambda_i(x^{(k)},\lambda_i^{(k)},r_i^{(k)})
  \quad\text{for}\quad i=0,1,\ldots,m.
\]
The following statements hold.
\begin{enumerate}
 \item \textup{Finite multiplier updates.}
 For any $i\in\{0,\ldots,m\}$, if $|\mathcal M_i|<+\infty$, then
 $\{\lambda_i^{(k)}\}$ is eventually constant and so it converges
 strongly. For $i=0$, the limit belongs to $\mathfrak A$.

 \item \textup{Finite penalty updates.}
 If $|\mathcal P_j|<+\infty$ for every $j=0,1,\ldots,m$, then, for each
 $i=0,1,\ldots,m$, the sequences $\{\lambda_i^{(k)}\}$ and
 $\{\widehat\lambda_i^{(k)}\}$ converge strongly to the same limit
 $\bar\lambda_i$.
\end{enumerate}
\end{theorem}
\begin{proof}
At each iteration, $\lambda_i^{(k+1)}$ is either $\lambda_i^{(k)}$ or
$\widehat\lambda_i^{(k)}$ when $k\in\mathcal M_i$. If
$|\mathcal M_i|<\infty$, the sequence $\{\lambda_i^{(k)}\}$ is eventually
constant, proving the first result.

We now prove the second result.  If $|\mathcal P_j|<\infty$ for every
$j=0,1,\ldots,m$, then {\tt update} is false after finitely many iterations
and there exists $\bar r_j\in\mathbb{R}_{>0}$ and $K\in\mathbb{N}$ for which
$r_j^{(k)}=\bar r_j$ for all $k\ge K$. The update rules then give
\[
  \|\widehat\lambda_i^{(k)}-\lambda_i^{(k)}\|_i
  \leq \bar r_i\tau_i^{(k)},\qquad
  \tau_i^{(k+1)}=q_i\tau_i^{(k)},\qquad
  q_i:=\min\left\{\frac{1}{\bar r_i},\theta_i\right\}^{\beta_i}\in(0,1).
\]
Thus $\tau_i^{(k)}=\tau_i^{(K)}q_i^{\,k-K}$ for $k\geq K$.
Whether the candidate is accepted or the stored multiplier is retained,
\[
 \|\lambda_i^{(k+1)}-\lambda_i^{(k)}\|_i
 \leq\|\widehat\lambda_i^{(k)}-\lambda_i^{(k)}\|_i
 \leq \bar r_i\tau_i^{(K)}q_i^{\,k-K}.
\]
Consequently, for integers $n>k\geq K$,
\[
 \begin{aligned}
 \|\lambda_i^{(n)}-\lambda_i^{(k)}\|_i
 &\leq\sum_{j=k}^{n-1}
       \|\lambda_i^{(j+1)}-\lambda_i^{(j)}\|_i
 &\leq \bar r_i\tau_i^{(K)}\sum_{j=k}^{n-1}q_i^{\,j-K}
 \leq\frac{\bar r_i\tau_i^{(K)}}{1-q_i}\,q_i^{\,k-K}.
 \end{aligned}
\]
The stored sequence is therefore Cauchy and converges strongly to some
$\bar\lambda_i$. Moreover,
\[
 \|\widehat\lambda_i^{(k)}-\bar\lambda_i\|_i
 \leq \bar r_i\tau_i^{(K)}q_i^{\,k-K}
      +\|\lambda_i^{(k)}-\bar\lambda_i\|_i
 \longrightarrow0,
\]
so the candidate sequence converges strongly to the same limit.
Furthermore, the initialization and risk update rule ensure that
$\lambda_0^{(k)}\in\mathfrak A$ for every $k$. Since $\mathfrak A$ is
closed, each risk limit asserted above belongs to $\mathfrak A$.
\end{proof}
\begin{corollary}\label{cor:finite.pen.asym.feas}
 Suppose Assumptions~(A0)-(A2) and (A5) hold and further suppose that
 the primal sequence $\{x^{(k)}\}$ is bounded.  If $|\mathcal P_j|<+\infty$
 for $j=0,1,\ldots,m$, then
 \[
   \lim_{k\to+\infty}\, d_C(Tx^{(k)}) = 0.
 \]
 In particular, every weak accumulation point $\bar x$ of $\{x^{(k)}\}$
 satisfies $T\bar x\in C$ and $T\bar{x}\in\partial I_{C_i}^*(\bar\lambda_i)$ for
  $i=1,\ldots,m$.
\end{corollary}
\begin{proof}
This results follows since
\[
  \begin{aligned}
  d_{C_i}(Tx^{(k)})
    &= \|Tx^{(k)}-\mathbf P_{C_i}(Tx^{(k)})\|_i \\
    &\le \|Tx^{(k)}-\mathbf P_{C_i}(Tx^{(k)}+(r^{(k)})^{-1}\lambda_i^{(k)})\|_i \\
    &= \frac{1}{r_i^{(k)}}\|\lambda_i^{(k)}-\widehat\lambda_i^{(k)}\|_i \le \tau_i^{(k)}
  \end{aligned}
\]
for $k$ sufficiently large.  Passing to the limit and noting that $r_i^{(k)}$
is constant for $k$ sufficiently large, we obtain
\[
  \lim_{k\to+\infty}\, d_{C_i}(Tx^{(k)}) = 0
\]
for $i=1,\ldots,m$.  The bounded regularity of $\{C_1,\ldots,C_m\}$ from
Assumption~(A5) ensures that $d_{C}(Tx^{(k)})\to 0$. The feasibility of weak
accumulation points then follows from the weak lower semicontinuity of
$x\mapsto d_{C}(Tx)$.  Finally, suppose $x^{(k_j)}\rightharpoonup\bar{x}$ as
$k_j\to+\infty$.  By the definition of $\widehat\lambda_i^{(k_j)}$, we have
that
\[
  r_i^{(k_j)}T x^{(k_j)} + \lambda_i^{(k_j)} \in \widehat\lambda_i^{(k_j)} + r_i^{(k_j)}\partial I_{C_i}^*(\widehat\lambda_i^{(k_j)})
  \quad\iff\quad
  T x^{(k_j)} + \frac{1}{r_i^{(k_j)}}(\lambda_i^{(k_j)}-\widehat\lambda_i^{(k_j)}) \in \partial I_{C_i}^*(\widehat\lambda_i^{(k_j)}).
\]
Since $\lambda_i^{(k_j)}\to \bar\lambda_i$,
$\widehat\lambda_i^{(k_j)}\to\bar\lambda_i$ and
$Tx^{(k_j)}\rightharpoonup T\bar x$,
we have that $T x^{(k_j)} + (r_i^{(k_j)})^{-1}(\lambda_i^{(k_j)}-\widehat\lambda_i^{(k_j)})\rightharpoonup T\bar x$
and therefore, $T\bar x\in\partial I_{C_i}^*(\bar\lambda_i)$ by
\cite[Proposition~16.36]{bauschke2017convex}.
\end{proof}

It is important to note that despite the strong statement in this theorem, we
are restricted to settings in which either the multipliers or the penalty
parameters are only updated \textit{finitely} many times. However, there
is no guarantee that either of these cases will appear in infinite dimensions.
We therefore present the following general asymptotic stationarity result,
which holds regardless of finite termination.

\begin{theorem}[Asymptotic Stationarity]
\label{thm:asymstat}
Let $\{(x^{(k)},\zeta^{(k)},\blam^{(k)},\mathbf{r}^{(k)})\}$ be an infinite
sequence generated by \Cref{algo:rales-method}, ignoring the stopping
conditions. Moreover, let Assumptions~(A0)-(A2) hold and assume the sequence of
iterates $\{x^{(k)}\}$ is bounded. Then the iterates asymptotically satisfy the
first-order optimality conditions in the sense that
\begin{equation}\label{eq:asymKKT}
\left\{\begin{aligned}
  &\lim_{k\to\infty} \|g(x^{(k)})\|_Y = 0 \\  
  &\limsup_{k\to\infty} \; \langle -(\mathbb{E}[\widehat{\lambda}_0^{(k)} f'(x^{(k)})] + g'(x^{(k)})^*\zeta^{(k)}), y - x^{(k)} \rangle_{X^*,X} \le 0 \quad \forall\, y \in T^{-1}(C).
\end{aligned}\right.
\end{equation}
Moreover, if $|\mathcal M_i \cup\mathcal P_i^c|=+\infty$ for some
$i=1,\ldots,m$, then
\[
  \liminf_{k\to\infty} d_{C_i}(Tx^{(k)}) = 0.
\]
\end{theorem}

\begin{proof}
By the subproblem stopping criterion of the algorithm, $g(x^{(k)})\to 0$ and
the derivatives of the augmented Lagrangian
\begin{equation}
  e^{(k)} := \mathbb{E}[\widehat{\lambda}_0^{(k)} f'(x^{(k)})] + \sum_{i=1}^m T^*\widehat{\lambda}_i^{(k)} + g'(x^{(k)})^*\zeta^{(k)} ,
\end{equation}
satisfy $\|e^{(k)}\|_{X^*} \le \varepsilon^{(k)} \to 0$, where
$\widehat{\lambda}_i^{(k)} := \Lambda_i(x^{(k)},\lambda_i^{(k)},r_i^{(k)})$
for $i=0,1,\ldots,m$. Rearranging this equality yields
\begin{equation}
  -(\mathbb{E}[\widehat{\lambda}_0^{(k)} f'(x^{(k)})] + g'(x^{(k)})^*\zeta^{(k)}) + e^{(k)} = \sum_{i=1}^m T^*\widehat{\lambda}_i^{(k)}.
\end{equation}
Applying both sides to $y-x^{(k)}$ for arbitrary $y\in T^{-1}(C)$, the result
follows from \cite[Prop.~5.7~\&~Rem.~5.8]{Antil2023}.
The final result, \cite[Lemma~5.1]{Antil2023} ensures that
$d_{C_i}(Tx^{(k)})\le\tau_i^{(k)}$ for all $k\in\mathcal P_i^c$ and therefore,
$d_{C_i}(Tx^{(k)})$ converges to zero if $|\mathcal P_i^c|=+\infty$.  On the
other hand, suppose $|\mathcal P_i|=+\infty$, then the multiplier update rule
in \Cref{algo:rales-update} ensures that
\[
  \frac{1}{r_i^{(k)}}\|\widehat\lambda_i^{(k)}\|_i
  \le \frac{\nu_i}{r_i^{(k)}}(r_i^{(k+1)})^{\gamma_i}
  \le \nu_i\eta_i^{\gamma_i}(r_i^{(k)})^{\gamma_i-1} \quad\forall\, k\in\mathcal M_i.
\]
Since $r_i^{(k)}\to+\infty$, \cite[Lemma~5.3]{Antil2023} proves that
$d_{C_i}(Tx^{(k)})$ converges to zero if $|\mathcal M_i|=+\infty$.
Combining these two results, proves the $d_{C_i}(Tx^{(k)})$ converges to zero
on $\mathcal M_i\cup\mathcal P_i^c$ when
$|\mathcal M_i\cup\mathcal P_i^c|=+\infty$ and the limit inferior result
follows.
\end{proof}

\Cref{eq:asymKKT} yields asymptotic stationarity under mild
assumptions. As observed in \cite{Antil2023}, we can strengthen this statement
by imposing stronger assumptions such as pseudomonotonicity, compactness, or
complete continuity of the derivative mappings to pass to the limit and obtain
a true KKT point.  Our final result addresses the remaining question about
feasibility that was partially answered in \Cref{cor:finite.pen.asym.feas}
and \Cref{thm:asymstat}.

\begin{theorem}[Feasibility of Weak Accumulation Points]\label{thm:weak.acc.feas}
Let the assumptions of \Cref{thm:asymstat} hold and suppose $\bar x$ is a weak
accumulation point of $\{x^{(k)}\}$ associated with the subsequence
$\{x^{(k_j)}\}$ and $\mathcal K=\{k_j\}\subseteq\mathbb{N}$. If
\begin{equation}\label{eq:inf.success}
  \left|\mathcal K \cap \bigcap_{i=1}^m (\mathcal M_i \cup \mathcal P_i^c)\right|=+\infty,
\end{equation}
then $T\bar x\in C$. On the other hand, if
Assumptions~(A4)-(A5) hold,
\begin{equation}\label{eq:weak.star.lagr}
  \theta^{(k_j)}:=\left(\sum_{i=1}^m r_i^{(k_j)}\right)^{-1} \to 0
  \quad\implies\quad \theta^{(k_j)}(\mathbb E[\Lambda_0(x^{(k_j)},\lambda_0^{(k_j)},r_0^{(k_j)})f'(x^{(k_j)})] + g'(x^{(k_j)})^*\zeta^{(k_j)})\rightharpoonup^* 0
\end{equation}
and $|\mathcal M_i \cup \mathcal P_i^c|<+\infty$ for at least one
$i=1,\ldots, m$, then there exists $\bar t_i\in(0,1]$ with
$\bar t_1 + \ldots + \bar t_m = 1$ 
and $\bar y_i\in C_i$ for
$i=1,\ldots,m$ such that
\begin{equation}\label{eq:xbar.conv.hull}
  T\bar x = \sum_{i=1}^m \bar t_i \bar y_i.
\end{equation}
\end{theorem}
\begin{proof}
Suppose that \eqref{eq:inf.success} holds. The proof of \Cref{thm:asymstat}
ensures that $d_{C_i}(Tx^{(k)})$ converges to zero on
$\mathcal M_i\cup\mathcal P_i^c$ and therefore, the weak lower semicontinuity
of $x\mapsto d_{C_i}(Tx)$ ensure that
\begin{equation}\label{eq:dist.liminf}
  0 = \liminf_{k_j\to+\infty}\,d_{C_i}(Tx^{(k_j)}) \ge d_{C_i}(T\bar x)
  \quad\iff\quad T\bar x\in C_i.
\end{equation}
According to \eqref{eq:inf.success}, this holds for $i=1,\ldots,m$ and so
$T\bar x\in C$.
Now, suppose that there exists $i$ for which $|\mathcal M_i\cup\mathcal P_i^c|<+\infty$.
In this case, we have that $|\mathcal M_i|<+\infty$ and
$|\mathcal P_i^c|<+\infty$ for some $i=1,\ldots,m$.  Consequently,
$r_\ell^{(k)}\to+\infty$ for all $\ell=1,\ldots,m$ and there exists
$K_1\in\mathbb{N}$ with $K_1>K_0$ such that
$\lambda_i^{(k+1)}=\lambda_i^{(k)}$ and
$r_\ell^{(k+1)}=\bar\eta r_\ell^{(k)}$ for all $k\ge K_1$ and
$\ell=1,\ldots,m$. Note that \eqref{eq:weak.star.lagr} yields
\[
  \begin{aligned}
  T^*\left(Tx^{(k_j)} + \sum_{\ell=1}^m \{\theta^{(k_j)}\lambda_\ell^{(k_j)}
  - t_\ell^{(k_j)}\mathbf P_{C_\ell}(z_\ell^{(k_j)})\}\right)
  &= \theta^{(k_j)}\sum_{\ell=1}^m T^*\widehat\lambda_\ell^{(k_j)} \\
  &= \theta^{(k_j)}(e^{(k_j)} - (\mathbb E[\widehat\lambda_0^{(k_j)}f'(x^{(k_j)})]+g'(x^{(k_j)})^*\zeta^{(k_j)})) \\ 
  &\rightharpoonup^* 0,
  \end{aligned}
\]
where $t_\ell^{(k_j)}=\theta^{(k_j)}r_\ell^{(k_j)}>0$, 
$t_1^{(k_j)}+\ldots+t_m^{(k_j)}=1$ for all $k_j$ and
$z_\ell^{(k_j)}=(r_\ell^{(k_j)})^{-1}\lambda_\ell^{(k_j)}+Tx^{(k_j)}$, 
and recall that $t_\ell^{(k_j)}=t_\ell^{(K_1)}=\bar t_\ell$ for all $k\ge K_1$.
Now, since $\{x^{(k_j)}\}$ weakly converges, it is bounded and so
$\{\mathbf P_{C_\ell}(z_\ell^{(k_j)})\}$ is also bounded for $\ell=1,\ldots,m$.
Consequently, there exists a further subsequence of
$\{\mathbf P_{C_\ell}(z^{(k_j)})\}$, which we do not relabel, that converges
weakly to some limit $\bar y_\ell\in C_\ell$ for $\ell=1,\ldots,m$.
Here, we have used the fact that $Z$ is Hilbert and $C_\ell$ is closed and convex
(hence, weakly closed).  Since  $r_\ell^{(k)}\to+\infty$,
\cite[Lemma~5.2]{Antil2023} yields $\theta^{(k_j)}\lambda_\ell^{(k_j)}\to 0$
for $\ell=1,\ldots,m$, which guarantees
\[
  \theta^{(k_j)}\sum_{\ell=1}^m T^*\widehat\lambda_\ell^{(k_j)}
  \rightharpoonup T^*\left(T\bar x - \sum_{\ell=1}^m \bar t_\ell\bar y_\ell\right) = 0.
\]
Finally, the injectivity of $T^*$ by Assumption~(A4) ensures that
$T\bar x = \sum_{\ell=1}^m \bar t_\ell\bar y_\ell$ as was to be shown.
\end{proof}
\begin{remark}
Since $\{x^{(k_j)}\}$ converges weakly to $\bar x$, it is bounded and
so \eqref{eq:weak.star.lagr} is satisfied if $\{\zeta^{(k_j)}\}$ is
bounded and $f'$ and $g'$ are bounded on bounded sets.  In addition,
the conclusion \eqref{eq:xbar.conv.hull} is a consequence of the
general Banach space setting and the independent penalty update in
\Cref{algo:rales-update}.  If, for example, $X$ is finite dimensional
or $m=1$,\ then \eqref{eq:xbar.conv.hull} implies that $T\bar x\in C$.
Alternatively, if $T$ is a compact operator, then again
\eqref{eq:xbar.conv.hull} implies that $T\bar x\in C$. In the case of
compact $T$ or finite-dimensional $X$, the result follows since
$Tx^{(k_j)}\to T\bar x$ and so the continuity of $\mathbf P_{C_i}$
ensures that $y_i=\mathbf P_{C_i}(T\bar x)$. We emphasize that
$\mathbf P_{C_i}$ is generally not weakly continuous and so we may only
have \eqref{eq:xbar.conv.hull}, i.e., $T\bar x$ is only guaranteed to
be contained in the convex hull of the union of $C_i$ for
$i=1,\ldots,m$, which can only happen if there exists at least one
$j=1,\ldots,m$ and $K\in\mathbb{N}$ sufficiently large such that 
\[
  \|\widehat\lambda_j^{(k)}-\lambda_j^{(k)}\|_j > r_j^{(k)}\tau_j^{(k)}
  \qquad\text{and}\qquad
  \|\widehat\lambda_j^{(k)}\|_j > \nu_j (r_j^{(k+1)})^{\gamma_j}
  \qquad\forall\, k\ge K.
\]
The second bound implies that $\|\widehat\lambda_j^{(k)}\|_j$
diverges at an exponential rate since $r_j^{(k+1)}=\bar\eta r_j^{(k)}$
for $k\ge K_0$.
\end{remark}

\subsection{Global Convergence under Exact Equality Constraint Feasibility}\label{ssec:conv-dual}
We now consider an idealized setting in which we assume we can compute
approximate minimizers and associated dual variables while maintaining strict
equality constraint feasibility, $\delta^{(k)} = 0$. To this end, we define
$X_{\rm ad} := \{x \in X \mid g(x) = 0\}$ to be the set of primal variables
satisfying the equality constraint.  In addition to strict equality constraint
feasibility, we strengthen the approximate stationarity check in line~2 of
\Cref{algo:rales-method} and assume instead that we can compute
$\epsilon$-minimizers of $L(\cdot,\blam^{(k)},\mathbf r^{(k)})$, which is often
only possible when $L(\cdot,\blam^{(k)},\mathbf r^{(k)})$ is convex.  Recall
that $x^{(k)}\in X_{\rm ad}$ is an $\epsilon$-minimizer of
$L(\cdot,\blam^{(k)},\mathbf r^{(k)})$ for $\epsilon\ge0$ if
\[
  L(x^{(k)},\blam^{(k)},\mathbf r^{(k)})-\epsilon \le \inf_{x\in X_{\rm ad}}\,L(x,\blam^{(k)},\mathbf r^{(k)}).
\]
The subsequent results generalize those in \cite{kouri2021primal}, with our next
result proving convergence of the primal variables in this idealized, yet
restrictive, setting.
For these results, we employ the following notation.
We equip $\mathfrak Z$ with its usual product Hilbert norm and, for a
positive penalty vector $\mathbf r$, use the weighted norm
\[
 \|\blam\|_{\mathfrak Z}^2
 :=\mathbb E[\lambda_0^2]+\sum_{i=1}^m\|\lambda_i\|_Z^2,
 \qquad
 \|\blam\|_{\mathbf r}^2
 :=\frac{1}{r_0}\mathbb E[\lambda_0^2]
   +\sum_{i=1}^m\frac{1}{r_i}\|\lambda_i\|_Z^2.
\]
For each iteration, set
\[
 r^{(k)}:=\min_{0\le i\le m}r_i^{(k)},\qquad
 R^{(k)}:=\max_{0\leq i\leq m}r_i^{(k)},\qquad\text{and}\qquad
 D^{(k)}\blam:=\left(\frac{\lambda_0}{r_0^{(k)}},\ldots,
                     \frac{\lambda_m}{r_m^{(k)}}\right).
\]
In particular,
$r^{(k)} \|\blam\|_{\mathbf r^{(k)}}^2\le\|\blam\|_{\mathfrak Z}^2 \leq R^{(k)}\|\blam\|_{\mathbf r^{(k)}}^2$
and $\|\blam\|_{\mathbf r^{(k)}}^2 = (D^{(k)}\blam,\blam)_{\mathfrak Z}$.
We also denote by $\bar{\mathbf r}\in(0,+\infty]^{m+1}$ the ``limit''
of $\{\mathbf r^{(k)}\}$, where we permit $\bar r_i = +\infty$.  Note
that the update scheme in \Cref{algo:rales-update} produces two cases:
either $\bar r_i<+\infty$ for all $i=0,1,\ldots,m$ or
$\bar r_i=+\infty$ for all $i=0,1,\ldots,m$.

\begin{theorem}[Primal Convergence]
\label{th:conv-primal}
Let $\{(x^{(k)}, \blam^{(k)}, \mathbf{r}^{(k)})\}$ be a sequence generated by
Algorithm~\ref{algo:rales-method} such that $x^{(k)} \in X_{\rm ad}$ satisfies:
\[
  L(x^{(k)}, \blam^{(k)}, \mathbf{r}^{(k)}) - \epsilon^{(k)}
  \le \inf_{x \in X_{\rm ad}} L(x, \blam^{(k)}, \mathbf{r}^{(k)}) 
\]
with $\epsilon^{(k)} \to \bar\epsilon\ge 0$.
Assume $X_{\rm ad}$ is weakly sequentially closed and that
Assumptions~(A0)-(A3) and (A5) hold. If the sequence $\{x^{(k)}\}$ admits a
weak accumulation point $\bar x$, then $\bar x$ is a $\bar \epsilon$-minimizer
of \eqref{eq:opt-prob}. That is, $\bar x\in X_{\rm ad}\cap T^{-1}(C)$ and
satisfies
\[
  \mathcal R(f(\bar x)) - \bar\epsilon \le \inf_{x\in X_{\rm ad}\cap T^{-1}(C)}\, \mathcal R(f(x)).
\]
In particular, if $\bar\epsilon=0$, then $\bar x$ is a minimizer of
$\mathcal R\circ f$ over $X_{\rm ad}\cap T^{-1}(C)$.
\end{theorem}

\begin{proof}
Let $\{x^{(k_j)}\}$ be a subsequence weakly converging to $\bar x$ in
$X_{\rm ad}$. We analyze the components of the augmented Lagrangian and
derive a bound to prove that $\bar x$ is a $\bar\epsilon$-minimizer.

\noindent \textit{(Risk Component $\Psi_0$):} 
Let $M>0$ denote the bound on the set $\mathfrak A$.
For $\mu \in \mathfrak{A}$, the dual variable update $\widehat\lambda_0^{(k)}$
satisfies
\[
  \mathbb E[(r_0^{(k)}f(x^{(k)})+\lambda_0^{(k)}-\widehat\lambda_0^{(k)})(\mu-\widehat\lambda_0^{(k)})] \le 0,
\]
which upon rearranging terms and using the boundedness of $\mathfrak A$ yields
\[
  r_0^{(k)}\mathbb E[f(x^{(k)})(\mu - \widehat\lambda_0^{(k)})]
  \le \mathbb E[(\widehat\lambda_0^{(k)}-\lambda_0^{(k)})(\mu-\widehat\lambda_0^{(k)})]
  \le 2M\|\widehat\lambda_0^{(k)}-\lambda_0^{(k)}\|_{L^2(\Omega,\mathcal{F},\mathbb{P})}.
\]
Passing to the supremum over $\mu\in\mathfrak A$ produces the bound
\[
  \mathcal R(f(x^{(k)})) - \mathbb E[f(x^{(k)})\widehat\lambda_0^{(k)}]
  \le \frac{2M}{r_0^{(k)}}\|\widehat\lambda_0^{(k)}-\lambda_0^{(k)}\|_{L^2(\Omega,\mathcal{F},\mathbb{P})}.
\]
Moreover, by the definition of $\Psi_0$, we have that
\[
  0 \le \mathcal R(f(x^{(k)})) - \Psi_0(x^{(k)},\lambda_0^{(k)},r_0^{(k)}).
\]
Combining the previous two bounds, and recalling that
$\frac{1}{2r_0^{(k)}}\|\lambda_0^{(k)}-\widehat\lambda_0^{(k)}\|_{L^2(\Omega,\mathcal{F},\mathbb{P})}^2 \le \frac{M}{r_0^{(k)}}\|\widehat\lambda_0^{(k)}-\lambda_0^{(k)}\|_{L^2(\Omega,\mathcal{F},\mathbb{P})}$,
yields the bounds
\begin{equation}\label{eq:bnd.psi0.risk}
  0 \le \mathcal R(f(x^{(k)})) - \Psi_0(x^{(k)},\lambda_0^{(k)},r_0^{(k)}) \le \frac{3M}{r_0^{(k)}}\|\widehat\lambda_0^{(k)}-\lambda_0^{(k)}\|_{L^2(\Omega,\mathcal{F},\mathbb{P})}.
\end{equation}

If the penalty parameters stabilize, then \Cref{th:conv-dual-finite} ensures
that $\|\widehat\lambda_0^{(k)}-\lambda_0^{(k)}\|_{L^2}\to0$ and so the
right-hand side of \eqref{eq:bnd.psi0.risk} converges to zero. Otherwise,
$r_0^{(k)}\to+\infty$ and the boundedness of $\mathfrak A$ ensures that the
right-hand side of \eqref{eq:bnd.psi0.risk} again converges to zero.
Consequently, we have that
\[
  \liminf_{k_j\to\infty} \, \Psi_0(x^{(k_j)},\lambda_0^{(k_j)},r_0^{(k_j)}) \ge \liminf_{k\to\infty}\, \mathcal R(f(x^{(k_j)})) \ge \mathcal R(f(\bar x)),
\]
where the final inequality follows from the weak lower semicontinuity of
$\mathcal R\circ f$.

\noindent \textit{(General Constraints $\Psi_i$):} For $i \in \{1, \dots, m\}$
and fixed $\mu \in Z$, the definition of $\Psi_i$ yields the lower bound:
\begin{align}\label{eq:psi_i_bnd}
  \Psi_i(x^{(k_j)}, \lambda_i^{(k_j)}, r_i^{(k_j)}) 
	\ge (\mu, Tx^{(k_j)})_Z - I_{C_i}^*(\mu) - \tfrac{1}{2r_i^{(k_j)}} \|\mu-\lambda_i^{(k_j)}\|_Z^2.
\end{align}
If $r_i^{(k_j)}\to+\infty$, then as argued in the proof of
\cite[Lemma~5.2]{Antil2023} and due to the fact that
$\gamma_i < 1/2$ for $i=1,\ldots,m$ (cf.\ \Cref{algo:rales-init}), the terms
$(r_i^{(k_j)})^{-1}\|\mu-\lambda_i^{(k_j)}\|_Z^2 \to 0$.
On the other hand, if $r_i^{(k_j)}\to\bar r_i<+\infty$, then
\Cref{th:conv-dual-finite} ensures that $\lambda_i^{(k_j)}\to\bar\lambda_i$,
$\widehat\lambda_i^{(k_j)}\to\bar\lambda_i$ and
$T\bar x\in\partial I_{C_i}^*(\bar\lambda_i)$ for $i=1,\ldots,m$.
Therefore, the right-hand side in \eqref{eq:psi_i_bnd} converges
to $(\mu,T\bar x)_Z-I_{C_i}^*(\mu)-\frac{1}{2\bar r_i}\|\mu-\bar\lambda_i\|_Z^2$.
For both cases, passing to the limit inferior as $k_j \to \infty$ in
\eqref{eq:psi_i_bnd} with $\mu$ fixed yields:
\[
  \liminf_{j \to \infty} \Psi_i(x^{(k_j)}, \lambda_i^{(k_j)}, r_i^{(k_j)}) \ge
  \left\{\begin{array}{ll}
  (\mu, T\bar x)_Z - I_{C_i}^*(\mu) & \text{if $r_i^{(k_j)}\to+\infty$} \\
  (\mu, T\bar x)_Z - I_{C_i}^*(\mu) - \frac{1}{2\bar r_i}\|\mu-\bar\lambda_i\|_Z^2 &\text{otherwise.}
  \end{array}\right.
\]
In the first case, taking the supremum over $\mu \in Z$ demonstrates
that the right-hand side is $I_{C_i}(T\bar x)$.
In the second case, we have that
\[
  \bar r_i T\bar x + \bar\lambda_i \in \bar\lambda_i + \bar r_i\partial I_{C_i}^*(\bar\lambda_i)
  \quad\iff\quad \bar\lambda_i=\Lambda_i(\bar x,\bar\lambda_i,\bar r_i)
\]
and therefore, \cite[Proposition~16.10]{bauschke2017convex} ensures that
$$\Psi_i(\bar x,\bar\lambda_i,\bar r_i)=(\bar\lambda_i,T\bar x)_Z -
I_{C_i}^*(\bar\lambda_i)=I_{C_i}(T\bar x).$$
Hence, in both cases, we have that
\[
  \liminf_{k_j \to \infty} \Psi_i(x^{(k_j)}, \lambda_i^{(k_j)}, r_i^{(k_j)}) \ge I_{C_i}(T\bar x).
\]

\noindent \textit{(Global Optimality):} 
Continuing, we observe that for any $x$
\[
\Psi_i(x,\lambda^{(k_j)}_i,r^{(k_j)}_i) \le \sup_{\mu \in Z}\{(\mu, Tx)_{Z} - I^*_{C_i}(\mu)\}.
\]
This, along with the observations above, provides the sequence of inequalities
\[
\aligned
\mathcal{R}(f(x)) + \sum_{i=1}^m \sup_{\mu \in Z}\{(\mu, Tx)_{Z} - I^*_{C_i}(\mu)\} 
&\ge L(x,\blam^{(k_j)},\mathbf{r}^{(k_j)}) \\
&\ge  L(x^{(k_j)},\blam^{(k_j)},\mathbf{r}^{(k_j)}) - \epsilon^{(k_j)}\\
&\ge \Psi_0(x^{(k_j)},\lambda_0^{(k_j)},r_0^{(k_j)})
  + \sum_{i=1}^m \Psi_i(x^{(k_j)},\lambda^{(k_j)}_i,r^{(k_j)}_i) - \epsilon^{(k_j)},
\endaligned
\]
for any $x\in X_{\rm ad}\cap T^{-1}(C)$.
Note that for such $x$, $\Psi_i(x, \lambda, r) \le 0$ holds. Then taking the
limit inferior over this inequality, in light of the weak lower semicontinuity
of $\mathcal{R} \circ f$, implies  
\[
  \mathcal{R}(f(x)) \ge \liminf_{j \to \infty} \left[ L(x^{(k_j)}, \blam^{(k_j)}, \mathbf{r}^{(k_j)}) - \epsilon^{(k_j)} \right]
                    \ge \mathcal{R}(f(\bar x)) + \sum_{i=1}^m I_{C_i}(T\bar x) - \bar\epsilon.
\]
The finiteness of the left-hand side implies $I_{C_i}(T\bar x) = 0$ for $i=1,\ldots,m$.
Hence, $T\bar x\in C$ and it follows that $\mathcal{R}(f(x))\ge\mathcal{R}(f(\bar x))-\bar\epsilon$
for all $x\in X_{\rm ad}\cap T^{-1}(C)$.
\end{proof}

Next, we turn our attention to the convergence of the dual variables. For this,
we define the \textbf{dual objective function} $v: \mathfrak{L} \to \mathbb R \cup \{\pm\infty\}$
by setting:
\[
  v(\boldsymbol{\lambda}) := \inf_{x \in X_{\rm ad}} \ell(x,\boldsymbol{\lambda}).
\]

\begin{proposition}[Dual Objective Function Properties]\label{prop:dual}
Under Assumptions~(A0) and (A2), the dual function $v$ is
concave and weakly upper semicontinuous. Moreover, if
\[
\mathfrak{S} := 
\Big\{ \blam \in \mathfrak{L} \;\Big|\; \operatorname*{arg\,min}_{x \in X_{\rm ad}} \ell(x,\boldsymbol{\lambda}) \ne \emptyset \Big\} \cap \Big\{ \boldsymbol{\lambda} \in \mathfrak{L} \; \Big| \; \partial I^*_{C_i}(\lambda_i) \neq\emptyset \; \forall i = 1,\dots, m\Big\} \neq \emptyset,
\]
then $(-v+I_{\mathfrak{L}}):\mathfrak{Z}\to\mathbb R\cup\{+\infty\}$ is
subdifferentiable on $\mathfrak S$ and $\partial(-v+I_{\mathfrak{L}})$ is
maximally monotone.
\end{proposition}

\begin{proof}
\noindent\textit{(Concavity).} Let $\blam \in \mathfrak{L}$ and recall that
\[
  v(\blam) = \inf_{x\in X_{\rm ad}}\,\ell(x,\blam) = \inf_{x\in X_{\rm ad}}\left\{ \mathbb E[\lambda_0 f(x)]+\sum_{i=1}^m(\lambda_i,Tx)_Z\right\} - \sum_{i=1}^m I_{C_i}^*(\lambda_i),
\]
where the first term on the right-hand side involves minimizing an affinely
parameterized objective function. The concavity of $v$ follows from this
and the convexity of $I_{C_i}^*$.

\noindent\textit{(Weak upper semicontinuity).} Let
$\{\blam^{(k)}\} \subset \mathfrak{L}$ such that
$\blam^{(k)} \rightharpoonup \blam$. Then we have for $v$:
\begin{align*}
v(\blam^{(k)}) &= \inf_{x \in X_{\rm ad}}\left\{\mathbb E[\lambda^{(k)}_0 f(x)] + \sum_{i=1}^m(\lambda^{(k)}_i, Tx)_{Z}\right\} - \sum_{i=1}^m I_{C_i}^*(\lambda_i^{(k)}) \\
  &\le \mathbb E[\lambda^{(k)}_0 f(x)] + \sum_{i=1}^m(\lambda^{(k)}_i, Tx)_{Z} - \sum_{i=1}^m I_{C_i}^*(\lambda_i^{(k)}) \quad \forall x \in X_{\rm ad}.
\end{align*}
Fixing an arbitrary $x \in X_{\rm ad}$, taking the limit superior over both
sides, and subsequently the infimum over $x \in X_{\rm ad}$ yields the
upper semicontinuity of $v$.

\noindent \textit{(Subdifferentiability).}
Let $\blam \in \mathfrak{S}$, $\bmu \in \mathfrak{L}$,
$x_{\blam} \in X_{\rm ad}$ be a minimizer of $\ell(\cdot,\blam)$,
and $\mathbf{y}=(y_1,\ldots,y_m)$ with $y_i\in\partial I_{C_i}^*(\lambda_i)$.
Define $V(\blam)=-v(\blam)+I_{\mathfrak{L}}(\blam)$.  The definition of $v$ yields
\begin{align*}
  V(\bmu) - V(\blam)
  &\ge -\ell(x_{\blam},\bmu)+\ell(x_{\blam},\blam) \\ 
  &\ge \mathbb E[(\mu_0 - \lambda_0)(-f(x_{\blam}))] + \sum_{i=1}^m(\mu_i - \lambda_i,-Tx_{\blam})_{Z}+\sum_{i=1}^m (y_i,\mu_i-\lambda_i)_Z.
\end{align*}
Note that this inequality is immediate if $\bmu\not\in\mathfrak L$
since $I_{\mathfrak{L}}(\bmu)=+\infty$. We may therefore conclude that
\[
  (-f(x_{\blam}), -Tx_{\blam}+y_1,\dots,-Tx_{\blam}+y_m) \in \partial V(\blam) \subset \mathfrak{Z}
\]
holds for any such $x_{\blam}$ and $\mathbf{y}$. Hence,
$V$ is subdifferentiable on the set $\mathfrak S$ and its
subdifferential is maximal monotone by
\cite[Th.~A]{RTRockafellar_1976a}.
\end{proof}

In the final result of this section, we relate the RALES algorithm,
\Cref{algo:rales-method} to the proximal point method applied to the
dual problem
\begin{equation}\label{eq:dual}
  \sup_{\blam\in\mathfrak L}\, v(\blam).
\end{equation}
In doing so, we are able to prove that the entire sequence of dual
variables $\{\blam^{(k)}\}$ converges weakly to a solution of \eqref{eq:dual}.
To achieve this result, we impose additional accuracy requirements on
the $\epsilon$-minimizers of $L(\cdot,\blam^{(k)},\mathbf r^{(k)})$.
As we will see, these requirements are relatively strict, but relax
the inexactness conditions from the classical analysis of augmented Lagrangian
methods \cite{RTRockafellar_1976a}.  In particular, we do not require squared
tolerances for the $\epsilon$-minimizers. One may be able to further relax
these requirements using modern developments in, e.g.,
\cite{eckstein2013practical}.

\begin{theorem}[Weak Dual Convergence]\label{th:conv-dual-ideal}
Let the assumptions of \Cref{prop:dual} hold and suppose that $\{x^{(k)}\}$
generated by \Cref{algo:rales-method} satisfies
\begin{equation}\label{eq:dual-primal-accuracy}
  L(x^{(k)},\blam^{(k)},\mathbf{r}^{(k)}) - \frac{\epsilon^{(k)}}{2 R^{(k)}}
  \le \inf_{x\in X_{\rm ad}} \; L(x,\blam^{(k)},\mathbf{r}^{(k)})
  \quad\text{with}\quad \sum_{k=1}^\infty \epsilon^{(k)}<\infty,\;\;\epsilon^{(k)}\ge 0.
\end{equation}
Moreover, define
$\widehat{\blam}^{(k)}:=(\widehat{\lambda}_0^{(k)},\widehat{\lambda}_1^{(k)},\ldots,\widehat{\lambda}_m^{(k)})$
and assume that for each $k$, the saddle-point equality holds:
\begin{equation}\label{eq:dual-minimax}
  \inf_{x \in X_{\rm ad}} L(x,\blam^{(k)},\mathbf r^{(k)})
  = \sup_{\bmu\in\mathfrak{L}} \left\{ v(\bmu) - \frac{1}{2} \|\bmu - \blam^{(k)}\|^2_{\mathbf{r}^{(k)}} \right\}.
\end{equation}
If $\{\widehat{\blam}^{(k)}\}$ is bounded, then both $\{\blam^{(k)}\}$ and
$\{\widehat{\blam}^{(k)}\}$ converge weakly to the same maximizer of the dual
problem.  In addition, if $|\mathcal{P}_j|<\infty$ for all $j=0,1,\ldots,m$,
then this convergence is strong.
\end{theorem}

\begin{proof}
Let $a_k:=\inf_{x\in X_{\rm ad}}L(x,\blam^{(k)},\mathbf r^{(k)})$ and note that
$a_k$ is finite by \eqref{eq:dual-primal-accuracy}.  For any
$\bar x\in X_{\rm ad}\cap T^{-1}(C)$, weak duality ensures that
$v(\bmu)\leq\mathcal R(f(\bar x))$ for every $\bmu\in\mathfrak L$. In
particular, $v-I_{\mathfrak L}$ is bounded above and its quadratic
regularization on the right-hand side of \eqref{eq:dual-minimax} is weakly
upper semicontinuous, coercive, and strongly concave.  Hence, the reflexivity
of $\mathfrak Z$ ensures that the maximization problem on the right-hand side
of \eqref{eq:dual-minimax}, has a unique maximizer in $\mathfrak L$, which we
denote by
\begin{equation}\label{eq:dual-prox-map}
 \mathbf p^{(k)}:=\operatorname*{arg\,max}_{\bmu\in\mathfrak L}
  \left\{v(\bmu)
       -\frac12\|\bmu-\blam^{(k)}\|_{\mathbf r^{(k)}}^2\right\}.
\end{equation}
Recall that $\mathbf p^{(k)}$ satisfies the optimality condition
\begin{equation}\label{eq:dual-prox-optimality}
  \mathbf w^{(k)}:=D_k(\blam^{(k)}-\mathbf p^{(k)})
    \in\partial (-v+I_{\mathfrak{L}})(\mathbf p^{(k)}).
\end{equation}

For $x\in X_{\rm ad}$ and $\bmu\in\mathfrak L$, set
\[
 H_k(x,\bmu):=\ell(x,\bmu)
             -\frac12\|\bmu-\blam^{(k)}\|_{\mathbf r^{(k)}}^2.
\]
By \Cref{prop:penalty_props}, $\widehat{\blam}^{(k)}$ is the unique
maximizer of $H_k(x^{(k)},\cdot)$, with maximum value
$L(x^{(k)},\blam^{(k)},\mathbf r^{(k)})$. Strong concavity with
respect to the $\|\cdot\|_{\mathbf r^{(k)}}$ norm yields
\[
 \frac12\|\widehat{\blam}^{(k)}-\mathbf p^{(k)}\|_{\mathbf r^{(k)}}^2
 \leq H_k(x^{(k)},\widehat{\blam}^{(k)})
      -H_k(x^{(k)},\mathbf p^{(k)}).
\]
Since $x^{(k)}\in X_{\rm ad}$, \eqref{eq:dual-minimax} gives
\[
 H_k(x^{(k)},\mathbf p^{(k)})
 \ge v(\mathbf p^{(k)})
    -\frac12\|\mathbf p^{(k)}-\blam^{(k)}\|_{\mathbf r^{(k)}}^2
 =a_k.
\]
Combining these inequalities with \eqref{eq:dual-primal-accuracy}, we obtain
\begin{equation}\label{eq:dual-prox-error}
 \|\mathbf p^{(k)}-\widehat{\blam}^{(k)}\|_{\mathbf r^{(k)}}^2
   \leq\frac{\epsilon^{(k)}}{R_k}\qquad\implies\qquad
 \|\mathbf p^{(k)}-\widehat{\blam}^{(k)}\|_{\mathfrak Z}
   \leq\sqrt{\epsilon^{(k)}}.
\end{equation}

We now consider two cases: $|\mathcal P_j|<\infty$ for every $j=0,1,\ldots,m$
and $|\mathcal P_j|=\infty$ for at least one $j=0,1,\ldots,m$.
In the first case, \Cref{th:conv-dual-finite} ensures the existence of
$\bar\blam$ for which $\blam^{(k)}\to\bar\blam$ and
$\widehat\blam^{(k)}\to\bar\blam$. Consequently, \eqref{eq:dual-prox-error}
yields $\mathbf p^{(k)}\to\bar\blam$ as well.
Since the penalty parameters are positive and nondecreasing,
$r^{(1)}>0$ satisfies $r_i^{(k)}\ge r^{(1)}$ for every $i$ and $k$. Hence,
\[
 \|\mathbf w^{(k)}\|_{\mathfrak Z}
 \leq \frac{1}{r^{(1)}}
 \|\boldsymbol{\lambda}^{(k)}-\mathbf p^{(k)}\|_{\mathfrak Z}
 \to 0
\]
and so $\mathbf w^{(k)}\to 0$.  Now, for fix $\mathbf z\in\mathfrak L$, the
optimality conditions \eqref{eq:dual-prox-optimality} and the definition of
the subdifferential ensure that
\begin{equation}\label{eq:dual-subgradient}
 v(\mathbf p^{(k)})\geq v(\mathbf z)
  +(\mathbf w^{(k)},\mathbf z-\mathbf p^{(k)})_{\mathfrak Z}.
\end{equation}
Passing to the limit and using upper semicontinuity of $v$ yield
\[
 v(\blam^*)\geq\limsup_{k\to\infty}v(\mathbf p^{(k)})\geq v(\mathbf z).
\]
Thus, $\blam^*$ maximizes $v$ over $\mathfrak L$.

To handle the second case, we recall that from iteration $K_0+1$ onward,
whenever the penalty parameters increase, all components are multiplied by the
same factor $\overline\eta>1$. Hence, for $k > K_0$,
\[
  r_i^{(k)}=\overline\eta^{n_k}r_i^{(K_0+1)} \qquad\text{for}\qquad i=0,\ldots,m,
\]
where $n_k$ counts the simultaneous penalty increases between the
$K_0+1$ and $k-1$ iterations with $n_{K_0+1}=0$. Since at least one
$\mathcal P_j$ is infinite, $n_k\to+\infty$, $r_i^{(k)}\to+\infty$ for
$i=0,\ldots,m$,  
\[
  \|\cdot\|_{\mathbf{r}^{(k)}}^2 = \overline\eta^{-n_k}\|\cdot\|_{\mathbf{r}^{(K_0+1)}}^2
  \qquad\text{and}\qquad 
  R^{(k)} = \overline{\eta}^{n_k} R^{(K_0+1)} \qquad\forall\, k > K_0. 
\]
Consequently, the boundedness of $\{\widehat\blam^{(k)}\}$ implies existence of
$K_1>K_0$ for which
\[
  \|\widehat\lambda_i^{(k)}\|_i \le \nu_i (r_i^{(k+1)})^{\gamma_i} \quad\forall\,i=0,1,\ldots,m,\;\; k\ge K_1
\]
and therefore, $\lambda_i^{(k+1)}=\widehat\lambda_i^{(k)}$
for all $i=0,1,\ldots,m$ and $k\ge K_1$.  This and \eqref{eq:dual-prox-error}
ensure that $\{\blam^{(k)}\}$ and $\{\mathbf p^{(k)}\}$ are bounded.
Moreover, this boundedness and the fact that $r_i^{(k)}\to+\infty$ for
$i=0,\ldots,m$ ensures that $\mathbf{w}^{(k)}\to 0$. Boundedness and the
fact that $\mathfrak{Z}$ is Hilbert guarantees the existence of a weakly
convergent subsequence $\widehat{\blam}^{(k_j)}\rightharpoonup\bar\blam$
and so $\mathbf p^{(k_j)}\rightharpoonup\bar\blam$ as well.  In addition,
$\bar\blam\in\mathfrak L$ since $\mathfrak L$ is weakly closed.
Applying the subgradient inequality \eqref{eq:dual-subgradient} along this
subsequence and using weak upper semicontinuity of $v$ ensures that, for every
$\mathbf z\in\mathfrak L$, we obtain
\[
 v(\bar\blam)
 \ge\limsup_{j\to\infty}v(\mathbf p^{(k_j)})
 \ge v(\mathbf z).
\]
Here we have used the by boundedness of $\mathbf p^{(k_j)}$ and 
the fact that $\mathbf w^{(k)}\to0$.  As a result, dual maximizers exist, and
every weak candidate accumulation point of $\{\widehat\blam^{(k)}\}$, and hence
$\{\blam^{(k)}\}$ and $\{\mathbf p^{(k)}\}$, is dual optimal.

Fix any dual solution $\mathbf z\in\mathfrak L$. The strong concavity of
$H_k(x^{(k)},\cdot)$ implies
\begin{align*}
 \frac12\|\widehat{\blam}^{(k)}-\mathbf z
             \|_{\mathbf r^{(k)}}^2
 &\le
 L(x^{(k)},\blam^{(k)},\mathbf r^{(k)})
 -\ell(x^{(k)},\mathbf z)
 +\frac12\|\mathbf z-\blam^{(k)}\|_{\mathbf r^{(k)}}^2\\
 &\le
 \frac{\epsilon^{(k)}}{2R^{(k)}}
 +\frac12\|\mathbf z-\blam^{(k)}\|_{\mathbf r^{(k)}}^2.
\end{align*}
The last inequality follows from
\eqref{eq:dual-primal-accuracy} and \eqref{eq:dual-minimax}; in particular,
$a_k \le v(\mathbf z) \le \ell(x^{(k)},\mathbf z)$.
Multiplying by $2\overline\eta^{n_k}$ and using the fact that
$\blam^{(k+1)}=\widehat{\blam}^{(k)}$ for $k\ge K_1$, we obtain
\[
 \|\blam^{(k+1)}-\mathbf z\|_{\mathbf r^{(K_0+1)}}^2
 \le
 \|\blam^{(k)}-\mathbf z\|_{\mathbf r^{(K_0+1)}}^2
 +\frac{\epsilon^{(k)}}{R^{(K_0+1)}}
 \qquad\forall\, k\ge K_1.
\]
From here, for all $k>k'\ge K_1$ we have that
\[
\begin{aligned}
  \|\blam^{(k)}-\mathbf z\|_{\mathbf r^{(K_0+1)}}^2 - \|\blam^{(k')}-\mathbf z\|_{\mathbf r^{(K_0+1)}}^2
  &= \sum_{\ell=k'}^{k-1} \{\|\blam^{(\ell+1)}-\mathbf z\|_{\mathbf r^{(K_0+1)}}^2-\|\blam^{(\ell)}-\mathbf z\|_{\mathbf r^{(K_0+1)}}^2\} \\
  &\le \frac{1}{R^{(K_0+1)}}\sum_{\ell=k'}^{k-1} \epsilon^{(\ell)}
   \le \frac{1}{R^{(K_0+1)}}\sum_{\ell=k'}^{+\infty} \epsilon^{(\ell)}.
\end{aligned}
\]
Since $\{\epsilon^{(k)}\}$ is summable we have that the tail sum
$\sum_{\ell=k'}^{+\infty}\epsilon^{(\ell)}\to0$ as $k'\to+\infty$.  Therefore,
$\{\|\blam^{(k)}-\mathbf z\|_{\mathbf r^{(K_0+1)}}^2\}$ is Cauchy and hence
converges, i.e.,
\[
  \lim_{k\to+\infty} \|\blam^{(k)}-\mathbf z\|_{\mathbf r^{(K_0+1)}}^2 = \kappa(\mathbf{z}) < +\infty,
\]
which holds for any dual solutions $\mathbf z$.  To finish the proof, we show
that all weak accumulation points of the sequence $\{\blam^{(k)}\}$ are the
same.  Suppose $\mathbf z_1$ and $\mathbf z_2$ are two distinct weak
accumulation points, which by the above arguments are dual optimal.  We
notice that
\[
 \|\blam^{(k)}-\mathbf z_2\|_{\mathbf r^{(K_0+1)}}^2
 -\|\blam^{(k)}-\mathbf z_1\|_{\mathbf r^{(K_0+1)}}^2
 =
 2(\blam^{(k)}-\mathbf z_1,\mathbf z_1-\mathbf z_2)_{\mathbf r^{(K_0+1)}}
 +\|\mathbf z_1-\mathbf z_2\|_{\mathbf r^{(K_0+1)}}^2,
\]
where the left-hand side has the limit
$\kappa(\mathbf z_2)-\kappa(\mathbf z_1)$ by the above arguments.
Consequently, the limit of $2(\blam^{(k)}-\mathbf z_1,\mathbf z_1-\mathbf z_2)_{\mathbf r^{(K_0+1)}}$
exists and is given by
\[
  2\lim_{k\to+\infty}(\blam^{(k)}-\mathbf z_1,\mathbf z_1-\mathbf z_2)_{\mathbf r^{(K_0+1)}} = \kappa(\mathbf z_2)-\kappa(\mathbf z_1)-\|\mathbf z_1-\mathbf z_2\|_{\mathbf r^{(K_0+1)}}^2.
\]
Since $\mathbf z_1$ is a weak accumulation point, the left-hand side must
be zero and so
\[
  \kappa(\mathbf z_2)-\kappa(\mathbf z_1) = \|\mathbf z_1-\mathbf z_2\|_{\mathbf r^{(K_0+1)}}^2 > 0.
\]
Interchanging the roles of $\mathbf z_1$ and $\mathbf z_2$ and applying the
above arguments gives
\[
  \kappa(\mathbf z_1)-\kappa(\mathbf z_2)
  =\|\mathbf z_1-\mathbf z_2\|_{\mathbf r^{(K_0+1)}}^2.
\]
Adding these identities yields
\[
  0=2\|\mathbf z_1-\mathbf z_2\|_{\mathbf r^{(K_0+1)}}^2,
\]
so $\mathbf z_1=\mathbf z_2$, contradicting their assumed
distinctness. Namely, all accumulation points are the same.
This further proves that
the entire sequence $\{\blam^{(k)}\}$ converges weakly to a dual
solution.
\end{proof}

\section{Applications}\label{sec:apps}

We consider the PDE-constrained optimization problem
\begin{subequations}\label{eq:semilinear}
\begin{align}
  &\min_{\substack{z\in L^2(D) \\ u_{\xi}\in H_0^1(D)\cap C_0(D)}}
     \mathcal{R}\left(\frac{1}{2}\int_D (u_{\xi}+1)^2\,\mathrm{d}x\right)
           + \frac{\alpha}{2}\int_D z^2\,\mathrm{d}x \\
  &\hspace{7mm} \operatorname{subject\;to}\qquad -10 \le z \le 10, \quad u_{\xi} \ge \psi \\
  &\hspace{2.7cm} \begin{aligned}
     -\Delta u_{\xi} + u_{\xi}^3 - \varrho u_{\xi} &= f_{\xi} + z &&\text{in $D$ a.s.} \\
     u_{\xi} &=0 && \text{on $\partial D$ a.s.}
  \end{aligned}\label{eq:semilinear-pde}
\end{align}
where $D=(0,1)^2$, $\alpha = 10^{-3}$, $\varrho = 1$, $\xi_i \sim \mathcal{N}(0,1)$ i.i.d.\ for $i=1,\ldots,200$, with
\begin{align}
  f_{\xi}(x) &= \sqrt{2} \sum_{i=1}^{100} \frac{\sin((i-\tfrac{1}{2})\pi x_1)}{(i-\tfrac{1}{2})\pi} \xi_{2i}
                                         + \frac{\sin((i-\tfrac{1}{2})\pi x_2)}{(i-\tfrac{1}{2})\pi} \xi_{2i-1} \quad{\text{and}}\\
  \psi(x) &= -\tfrac{2}{3} + \tfrac{1}{2} \min\{x_1+x_2, \min\{1+x_1-x_2, \min\{1-x_1+x_2, 2-x_1-x_2\}\}\} .
\end{align}
\end{subequations}
To analyze the PDE \eqref{eq:semilinear-pde}, we employ the
parametrized operator $A_t\in\mathcal{L}(H_0^1(D),H^{-1}(D))$ defined by
\[
  \langle A_t u, v\rangle_{H^{-1}(D),H_0^1(D)} := \int_D \{\nabla u\cdot\nabla v-tuv\}\,\mathrm{d}x
  \qquad\text{for}\qquad t\in\mathbb{R}.
\]
Note that $A_{\varrho}$ is coercive on $H_0^1(D)$ since
\[
  \lambda_{\min}\|u\|_{L^2(D)}^2 \le \|u\|_{H^1_0(D)}^2 \quad\forall\, u\in H^1_0(D),
\]
where $\lambda_{\min}=2\pi^2$ is the minimum eigenvalue of $A_0$.  In
particular,
\[
  \langle A_\varrho u,u\rangle_{H^{-1}(D),H_0^1(D)} = \|u\|_{H^1_0(D)}^2 - \varrho\|u\|_{L^2(D)}^2
  \ge (1-\lambda_{\min}^{-1}\varrho)\|u\|_{H^1_0(D)}^2 \quad\forall\, u\in H^1_0(D).
\]
Since $\varrho=1$, the constant in the lower bound is positive and $A_\varrho$
is coercive. Additionally, recall that $H^1_0(D)$ is embedded in $L^6(D)$ and
so $u^3\in L^2(D)$ for $u\in H_0^1(D)$ and $N(u):=u^3$ is a monotone
nonlinearity.  Therefore, \cite[Theorem~2.1]{kouri2020risk} ensures that
solutions $u_\xi$ to \eqref{eq:semilinear-pde} exist, are unique and satisfy
$u_\xi\in L^p(\Omega,\mathcal{F},\mathbb{P};H^1_0(D))$ for any $1\le p<\infty$.
The integrability here follows since
$\|f_\xi\|_{L^2(D)}\in L^p(\Omega,\mathcal{F},\mathbb{P})$ for any
$1\le p<\infty$. Finally, fix $\xi\in\real^{200}$. As previously discussed, we
have that $u_\xi^3\in L^2(D)$ and so
$A_0u_\xi=f_\xi+z-u_\xi^3+u_\xi\in L^2(D)$.
Consequently, $u_\xi\in H^2(D)\cap C_0(D)$ for each $\xi\in\real^{200}$.

The spaces and maps corresponding to \eqref{eq:semilinear} are
listed below. Feasibility is addressed separately in
\Cref{rem:sample-feasibility}.
\begin{equation*}
\begin{aligned}
x&=(u_\xi, z) \in X=L^4(\Omega,\mathcal{F},\mathbb{P};H^{1}_0(D)\cap C_0(D))\times L^2(D), \\[5pt]
g(x)&=Au_\xi+N(u_\xi)-f_\xi-z \in Y=L^{4/3}(\Omega,\mathcal{F},\mathbb{P};H^{-1}(D)), \\[5pt]
Tx&=(u_\xi,z) \in Z=L^2(\Omega,\mathcal{F},\mathbb{P};L^2(D))\times L^2(D), \\[5pt]
C_1&=\{u\in L^2(\Omega,\mathcal{F},\mathbb{P};L^2(D))\,\vert\,u\ge\psi\text{ a.s.}\}\times L^2(D), \\[5pt]
C_2&=L^2(\Omega,\mathcal F,\mathbb P;L^2(D))\times \{z\in L^2(D)\,\vert\,-10\le z\text{ a.e.}\}, \\[5pt]
C_3&=L^2(\Omega,\mathcal F,\mathbb P;L^2(D))\times \{z\in L^2(D)\,\vert\,z\le 10\text{ a.e.}\} . \\[5pt]
\end{aligned}
\end{equation*}

\begin{remark}[Feasibility of \eqref{eq:semilinear}]
\label{rem:sample-feasibility}
  Verifying the feasibility of almost sure state constraints is nontrivial
  and often not possible.  This theoretical challenge can lead to unknowingly
  formulating an infeasible problem, potentially causing problems with its
  numerical solution.  In fact, \eqref{eq:semilinear} is infeasible. To see
  this, consider the first eigen function of $A_0$,
  $\phi(x)=\sin(\pi x_1)\sin(\pi x_2)$.  Testing \eqref{eq:semilinear-pde}
  with $\phi$ for fixed $\xi\in\mathbb{R}^{200}$ yields
  \[
    \begin{aligned} 
    (f_\xi, \phi)_{L^2(D)} &= \langle A_\varrho u_\xi, \phi\rangle_{H^{-1}(D),H_0^1(D)} + (u_\xi^3-z,\phi)_{L^2(D)} \\ 
      &= (u_\xi^3+(2\pi^2-\varrho)u_\xi-z,\phi)_{L^2(D)} =: (g(u_\xi) - z,\phi)_{L^2(D)}.
    \end{aligned} 
  \]
  
  Note that the polynomial $g(t)=t^3+(2\pi^2-\varrho)t$ is increasing and
  $\psi\ge -2/3$ a.e.  Consequently, if $u_\xi\ge\psi$ a.e.\ and $z\le 10$
  a.e., then
  \[
    (f_\xi, \phi)_{L^2(D)} \ge (g(-2/3)-10)\int_D \phi\,\textup{d}x =: M.
  \]
  Since this holds for all $\xi\in\mathbb{R}^{200}$, we conclude that
  \[
    \mathbb P\left((f_\xi,\phi)_{L^2(D)} < M\right) = 0,
  \]
  which cannot happen since $\xi_i\sim\mathcal N(0,1)$ can take values with
  arbitrarily large magnitude.  This contradiction demonstrates that
  \eqref{eq:semilinear} is infeasible.  
  
  For a given finite sample set, a sufficient condition for feasibility
at the continuous PDE level is $f_\xi\ge-9$ on $D$ at every sample.
Choose the common control $z=9$. Testing each state equation with
$v=u_\xi^-:=\max\{-u_\xi,0\} \ge 0$ and using Poincar\'e's inequality gives
\[
 \left(1-\frac{\varrho}{2\pi^2}\right)
 \|\nabla v\|_{L^2(D)}^2+\|v\|_{L^4(D)}^4\le0.
\]
Since $\varrho = 1$ and  $2\pi^2 \ge 1$, we obtain that $v = 0$. 
Moreover, since $\psi\le-1/6$, it follows that $u_\xi\ge0>\psi$, so this
control is feasible at those samples.
\end{remark}

We study two risk measures: the mean plus semideviation (MPSD)
\begin{align}
  \mathcal{R}[F] = \mathbb{E}[F] + p^{\text{MP}} \mathbb{E}[\max\{0,F-\mathbb{E}[F]\}]
\end{align}
where $p^{\text{MP}}\in[0,1)$ and the convex combination of the mean and the
Conditional Value-at-Risk (CVaR)
\begin{align}
  \mathcal{R}[F] = (1-p^{\text{CV}}) \mathbb{E}[F] + p^{\text{CV}} \text{CVaR}_\beta(F),
\end{align}
where $p^{\text{CV}} \in [0,1]$ and $\beta\in[0,1)$ is the confidence level.
By defining the random-valued objective function
\[
f(x) = \frac{1}{2}\int_D (u_{\xi}+1)^2\,\mathrm{d}x + \frac{\alpha}{2}\int_D z^2\,\mathrm{d}x
\]
the problem matches our theoretical form~\eqref{eq:opt-prob}, because the risk measures
are translation equivariant (i.e., $\mathcal{R}(X + c) = \mathcal{R}(X) + c$ for all
$X\in L^2(\Omega,\mathcal{F},\mathbb{P})$ and $c\in\real$) and the control penalty
$\frac{\alpha}{2}\int_D z^2\,\mathrm{d}x$ is deterministic.

\begin{figure}[hbt!]
  \centering
  \includegraphics[width=0.34\textwidth]{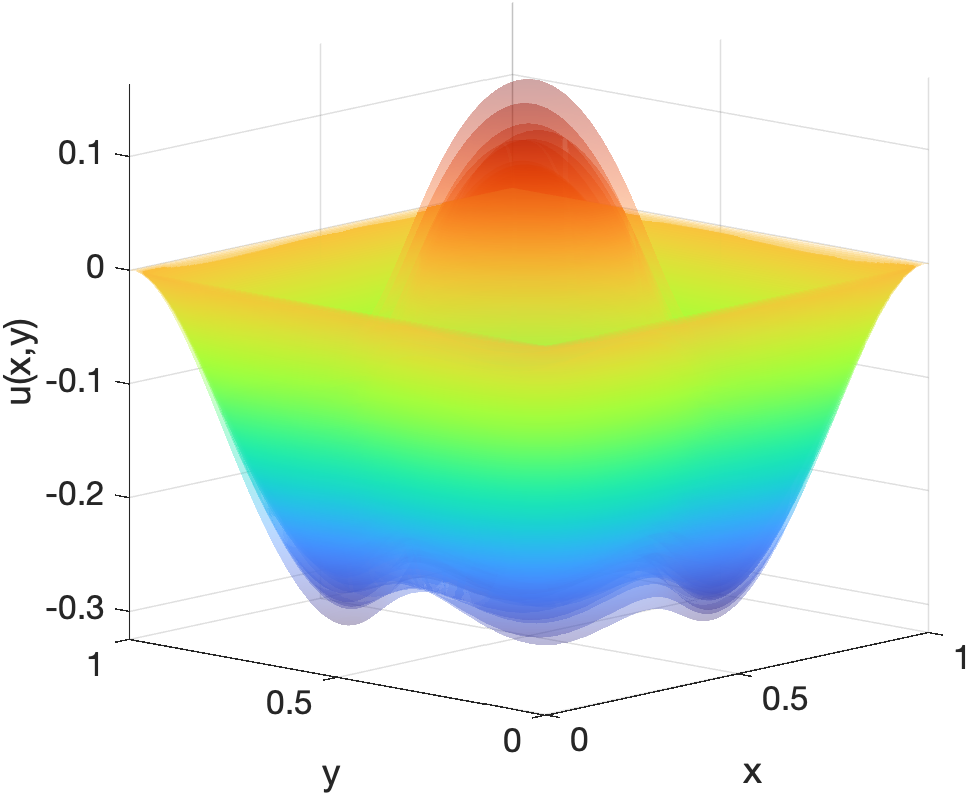}
  \includegraphics[width=0.30\textwidth]{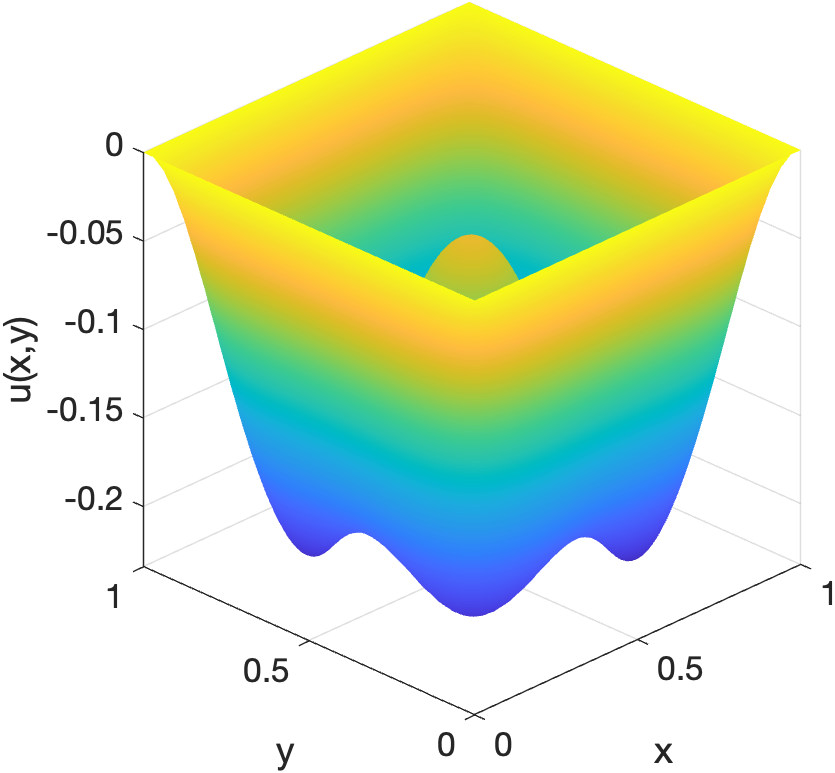}
  \includegraphics[width=0.32\textwidth]{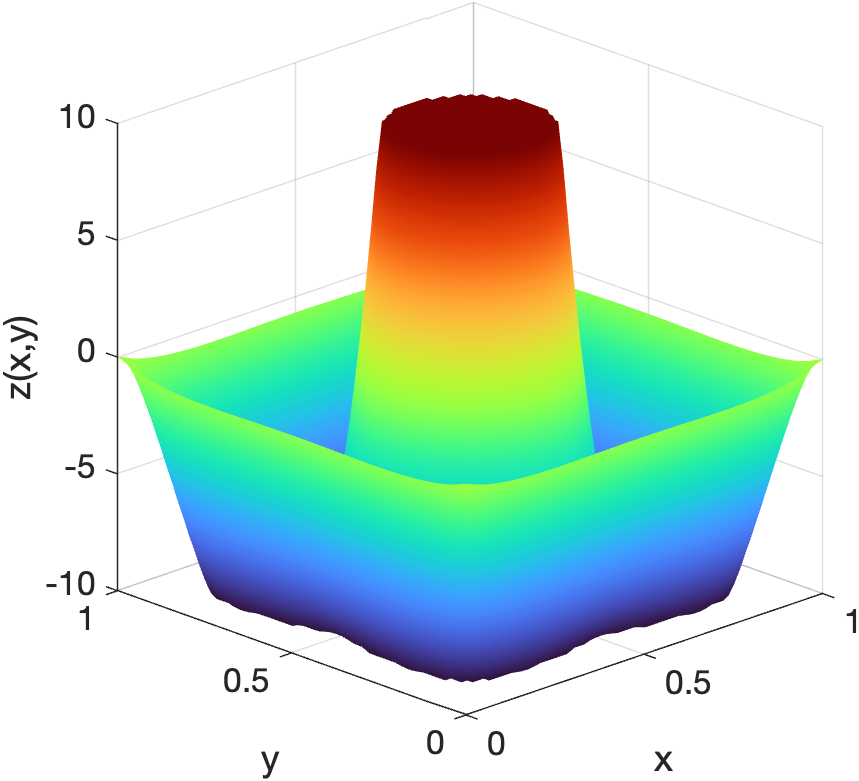}
\caption{Optimal state ``cloud'' for all samples (left), its mean (middle),
         and the optimal control (right).}
\label{fig:solution}
\end{figure}

To numerically solve \eqref{eq:semilinear}, we discretize in space using
continuous piecewise linear (P1) finite elements and we approximate the
risk measures using finitely many samples via Sample Average Approximation
(SAA). Figure~\ref{fig:solution} illustrates the solution of the optimal
control problem for a 64$\times$64 mesh and 100 samples of the normal
distribution.

We evaluate the numerical performance of RALES based on algorithmic and accuracy criteria.
Our algorithmic criteria include the numbers of outer RALES iterations, labeled \texttt{AL},
the numbers of SQP iterations, labeled \texttt{SQP}, and the numbers of projected conjugate
gradient iterations, labeled \texttt{CG}, needed for the solution of the tangential subproblems,
which represent the computational bulk of SQP and therefore RALES.
All reported iteration numbers are \emph{totals}.
Our accuracy criteria include the equality constraint norm,
$$\texttt{normg}=\|g(x^{(k)})\|_{L^{2}(\Omega;H^{-1}(D))},$$
the norm of the gradient of the subproblem Lagrangian,
$$\texttt{grad-lag}=\|L_x'(x^{(k)},\blam^{(k)},\mathbf{r}^{(k)}) + g'(x^{(k)})^*\zeta^{(k)}\|_{X^*},$$
the measure of inequality constraint infeasibility,
$$\texttt{feas}=\max_i d_{C_i}(Tx^{(k)}),$$
and the dual risk metric,
$$\texttt{dual-risk}=\frac{1}{r_0^{(k)}}\|\lambda_0^{(k)}-\Lambda_0(x^{(k)},\lambda_0^{(k)},r_0^{(k)})\|_{L^2(\Omega,\mathcal{F},\mathbb{P})}.$$
Throughout, we set the tolerances $\varepsilon_*=\delta_*=\tau_*=\sigma_*=10^{-6}$.
The multiplier stopping test in \Cref{algo:rales-method} uses
\[
 \texttt{dual-mult}:=\max_{1\le i\le m}
 \frac{1}{r_i^{(k)}}\|\Lambda_i(x^{(k)},\lambda_i^{(k)},r_i^{(k)})
              -\lambda_i^{(k)}\|_i.
\]
The projection estimate above gives $\texttt{feas}\le\texttt{dual-mult}$.
Thus a small value of \texttt{feas} alone does not verify this stopping
test.

First, we investigate the iteration scaling as the computational mesh is refined.
The number of samples of the normal distribution is set to 100.
Tables~\ref{tab:MPSD-mesh} and~\ref{tab:CVAR-mesh} show that the \texttt{AL}
iteration numbers are nearly mesh independent and that the \texttt{SQP} and
\texttt{CG} iteration numbers increase moderately with mesh refinement.
Additionally, high accuracy is achieved at the final RALES iterate, in all reported
measures of progress.
This observation is consistent with the ALESQP
result~\cite[Table 4, tol 1e-6]{Antil2023}, with the additional observation
that the measure of dual risk is near machine zero (and often equals zero).
We note that the linear (augmented) systems in the tangential step computation and
other trust-region composite-step SQP subroutines are solved using flexible GMRES
preconditioned by an approximate Schur complement~\cite{rees2010optimal} based on
the state Jacobian matrix and its transpose.
The state Jacobian matrix is block-diagonal, i.e., the inverse of the approximate
Schur complement can be applied efficiently in parallel across all samples,
see~\cite[Fig.~1]{kouri2018inexact}.
We also note that, in addition to the iterative (inexact) solution of augmented
systems, the applications of the state Jacobian inverses and their transposes
within the augmented system preconditioner need not be very accurate.

\begin{table}[htb!]
\centering
\begin{tabular}{ c | c c c c c c c }
\hline
 {\tt mesh} & {\tt AL} & {\tt SQP} & {\tt CG} & {\tt normg} & {\tt grad-lag} & {\tt feas} & {\tt dual-risk} \\
\hline
 {\tt   32$\times$32} & {\tt 17} &  {\tt 85} & {\tt 200} & {\tt 2.87e-10} & {\tt 2.72e-07} & {\tt 8.83e-09} & {\tt 0.00e+00} \\
 {\tt   64$\times$64} & {\tt 18} &  {\tt 94} & {\tt 266} & {\tt 1.42e-10} & {\tt 7.92e-08} & {\tt 5.92e-08} & {\tt 0.00e+00} \\
 {\tt 128$\times$128} & {\tt 21} & {\tt 103} & {\tt 283} & {\tt 7.51e-11} & {\tt 2.06e-08} & {\tt 6.44e-08} & {\tt 0.00e+00} \\
 {\tt 256$\times$256} & {\tt 24} & {\tt 175} & {\tt 433} & {\tt 2.19e-13} & {\tt 2.27e-09} & {\tt 9.81e-07} & {\tt 0.00e+00} \\
\hline
\end{tabular}
\caption{Mean plus semideviation: Iteration numbers and feasibility measures at the final iterate.
         The risk parameter is set to $p^{\text{MP}} = 0.8$.}
\label{tab:MPSD-mesh}
\end{table}

\begin{table}[htb!]
\centering
\begin{tabular}{ c | c c c c c c c }
\hline
 {\tt mesh} & {\tt AL} & {\tt SQP} & {\tt CG} & {\tt normg} & {\tt grad-lag} & {\tt feas} & {\tt dual-risk} \\
\hline
 {\tt   32$\times$32} & {\tt 18} &  {\tt 77} & {\tt 215} & {\tt 2.38e-11} & {\tt 3.10e-08} & {\tt 5.32e-09} & {\tt 0.00e+00} \\
 {\tt   64$\times$64} & {\tt 18} & {\tt  98} & {\tt 263} & {\tt 2.79e-10} & {\tt 3.98e-07} & {\tt 7.01e-08} & {\tt 3.31e-10} \\
 {\tt 128$\times$128} & {\tt 21} & {\tt 133} & {\tt 351} & {\tt 4.71e-11} & {\tt 5.44e-07} & {\tt 8.12e-08} & {\tt 0.00e+00} \\
 {\tt 256$\times$256} & {\tt 24} & {\tt 160} & {\tt 425} & {\tt 2.48e-11} & {\tt 9.02e-08} & {\tt 8.89e-07} & {\tt 0.00e+00} \\
\hline
\end{tabular}
\caption{Mean plus CVaR: Iteration numbers and feasibility measures at the final iterate.
         The risk parameter is set to $p^{\text{CV}} = 0.75$ and the confidence level is $\beta = 0.8$.}
\label{tab:CVAR-mesh}
\end{table}

Second, we investigate the iteration scaling as the number of samples of the normal
distribution increases.
The mesh resolution is set to~64$\times$64.
Tables~\ref{tab:MPSD-sample} and~\ref{tab:CVAR-sample} show that, again, the \texttt{AL}
iteration numbers are nearly sample-size independent, the \texttt{SQP} and \texttt{CG}
iteration numbers are mildly sample-size dependent, and that RALES meets stringent
accuracy criteria. 

\begin{table}[htb!]
\centering
\begin{tabular}{ c | c c c c c c c }
\hline
 {\tt samples} & {\tt AL} & {\tt SQP} & {\tt CG} & {\tt normg} & {\tt grad-lag} & {\tt feas} & {\tt dual-risk} \\
\hline
 {\tt 100} & {\tt 18} &  {\tt 94} & {\tt 266} & {\tt 1.42e-10} & {\tt 7.92e-08} & {\tt 5.92e-08} & {\tt 0.00e+00} \\
 {\tt 200} & {\tt 17} & {\tt 139} & {\tt 310} & {\tt 6.77e-11} & {\tt 5.61e-08} & {\tt 8.22e-07} & {\tt 0.00e+00} \\
 {\tt 400} & {\tt 17} & {\tt 148} & {\tt 324} & {\tt 1.60e-10} & {\tt 3.92e-07} & {\tt 7.67e-07} & {\tt 4.46e-10} \\
 {\tt 800} & {\tt 22} & {\tt 152} & {\tt 443} & {\tt 6.27e-12} & {\tt 9.18e-08} & {\tt 8.48e-08} & {\tt 3.09e-12} \\
\hline
\end{tabular}
\caption{Mean plus semideviation: Iteration numbers and feasibility measures at the final iterate.
         The risk parameter is set to $p^{\text{MP}} = 0.8$.}
\label{tab:MPSD-sample}
\end{table}

\begin{table}[htb!]
\centering
\begin{tabular}{ c | c c c c c c c }
\hline
 {\tt samples} & {\tt AL} & {\tt SQP} & {\tt CG} & {\tt normg} & {\tt grad-lag} & {\tt feas} & {\tt dual-risk} \\
\hline
 {\tt 100} & {\tt 18} & {\tt  98} & {\tt 263} & {\tt 2.79e-10} & {\tt 3.98e-07} & {\tt 7.01e-08} & {\tt 3.31e-10} \\
 {\tt 200} & {\tt 18} & {\tt 107} & {\tt 280} & {\tt 6.07e-11} & {\tt 2.83e-07} & {\tt 1.23e-08} & {\tt 0.00e+00} \\
 {\tt 400} & {\tt 18} & {\tt 123} & {\tt 319} & {\tt 1.76e-11} & {\tt 2.66e-07} & {\tt 1.66e-07} & {\tt 1.14e-09} \\
 {\tt 800} & {\tt 22} & {\tt 145} & {\tt 391} & {\tt 2.55e-13} & {\tt 5.24e-08} & {\tt 1.47e-07} & {\tt 8.22e-12} \\
\hline
\end{tabular}
\caption{Mean plus CVaR: Iteration numbers and feasibility measures at the final iterate.
         The risk parameter is set to $p^{\text{CV}} = 0.75$ and the confidence level is $\beta = 0.8$.}
\label{tab:CVAR-sample}
\end{table}

\begin{figure}[hbt!]
  \centering
  \includegraphics[width=0.6\textwidth]{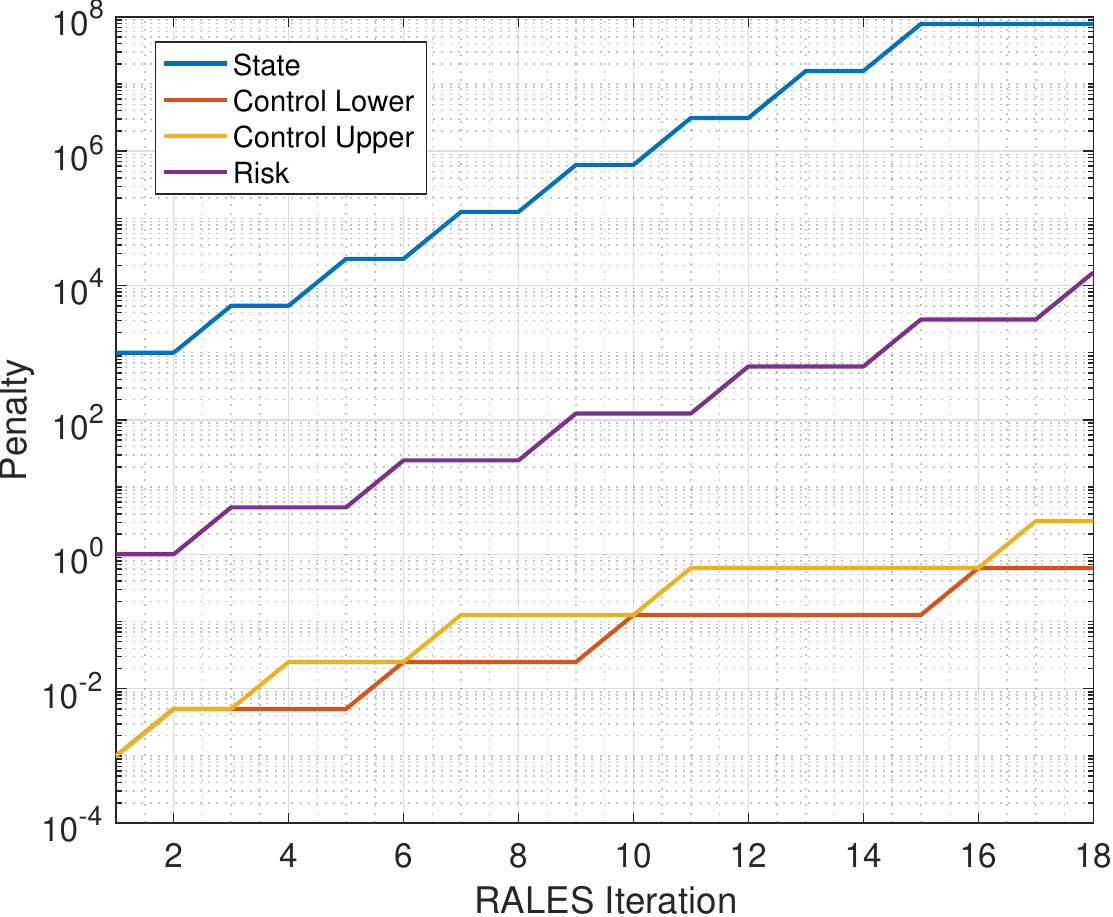}
\caption{Augmented Lagrangian penalty parameters across RALES iterations.}
\label{fig:penalties}
\end{figure}

Third, for the case of mean plus CVaR, in Figure~\ref{fig:penalties}
we report the evolution of the augmented Lagrangian penalty parameters
across RALES iterations, corresponding to Table~\ref{tab:CVAR-mesh},
row 2 (64$\times$64 mesh, 100 samples).
We note that while the state constraint penalty reaches large values,
near $10^8$, the risk penalty exhibits slower growth, and the control
penalties remain relatively modest (less than 10).

\begin{figure}[htb!]
     \centering
     \begin{subfigure}[b]{0.45\textwidth}
         \centering
         \includegraphics[width=0.8\textwidth]{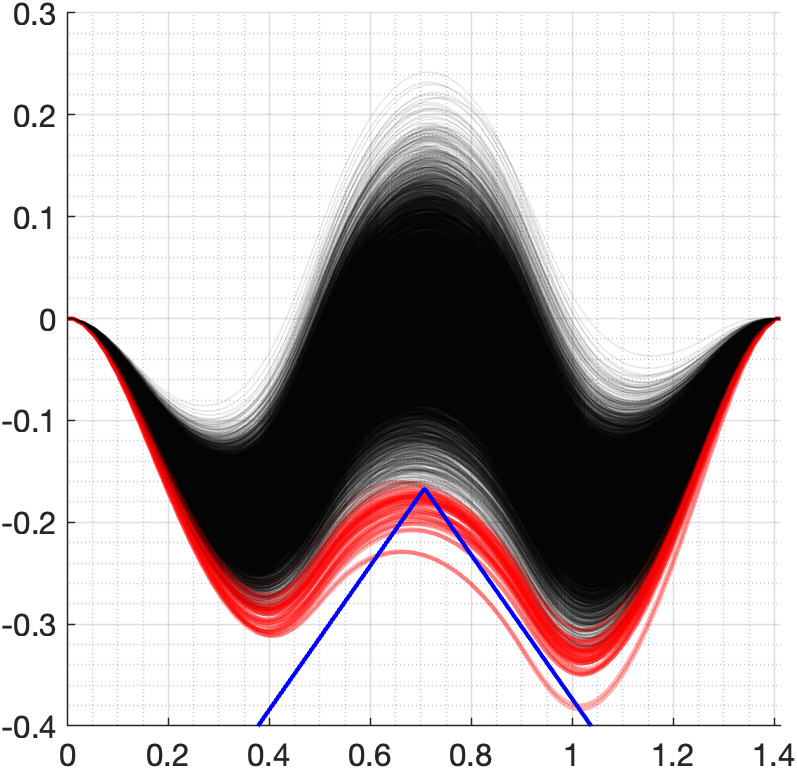} 
         \caption{100 original samples: 40 state constraint violations
                  among 10,000 new samples.}
     \end{subfigure}
     \hfill
     \begin{subfigure}[b]{0.45\textwidth}
         \centering
         \includegraphics[width=0.8\textwidth]{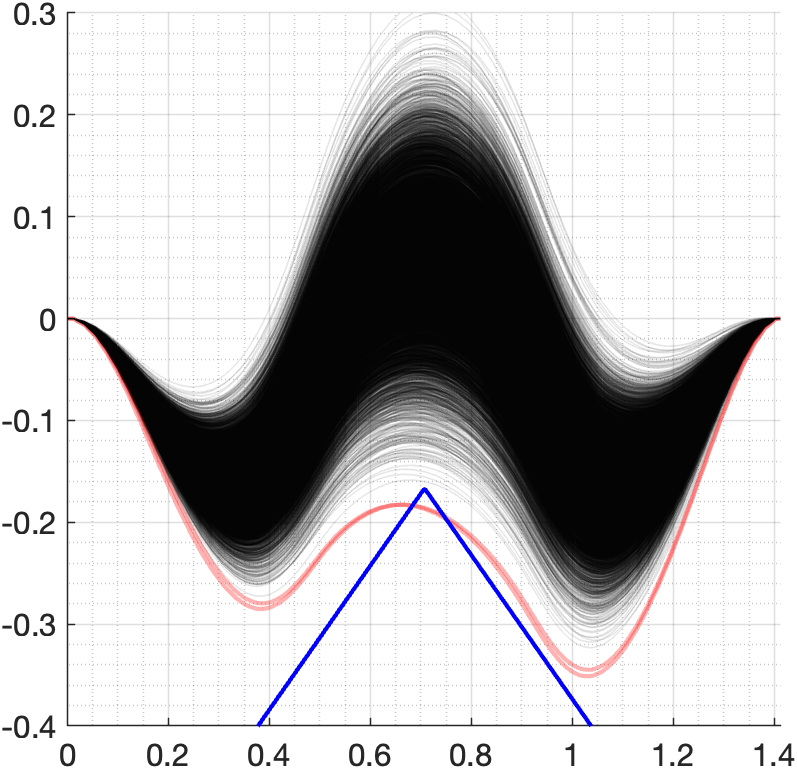} 
         \caption{200 original samples: 2 state constraint violations
                  among 10,000 new samples.}
     \end{subfigure}

     \vspace{1em} 

     \begin{subfigure}[b]{0.45\textwidth}
         \centering
         \includegraphics[width=0.8\textwidth]{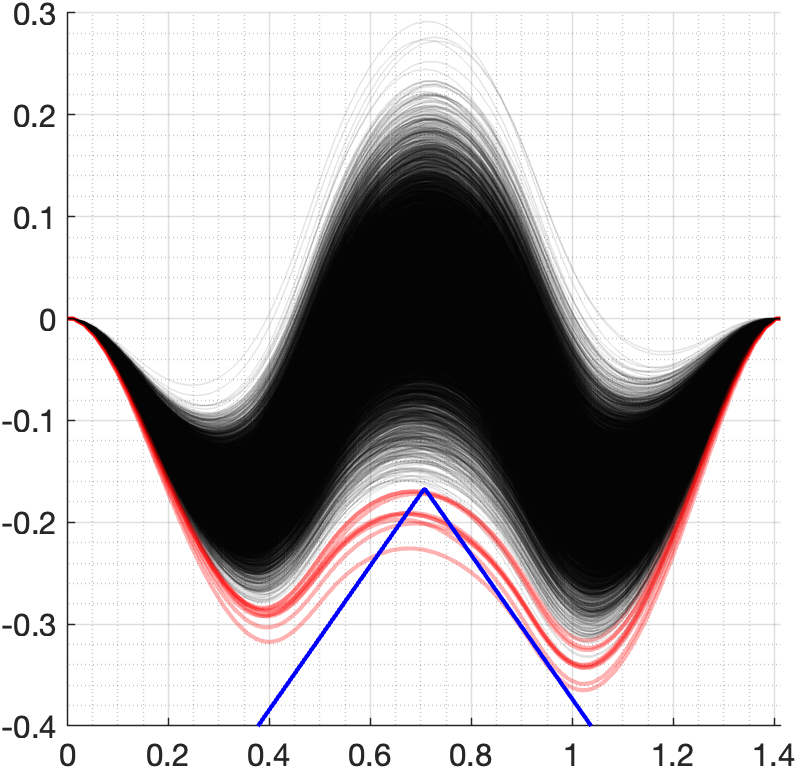} 
         \caption{400 original samples: 9 state constraint violations
                  among 10,000 new samples.}
     \end{subfigure}
     \hfill
     \begin{subfigure}[b]{0.45\textwidth}
         \centering
         \includegraphics[width=0.8\textwidth]{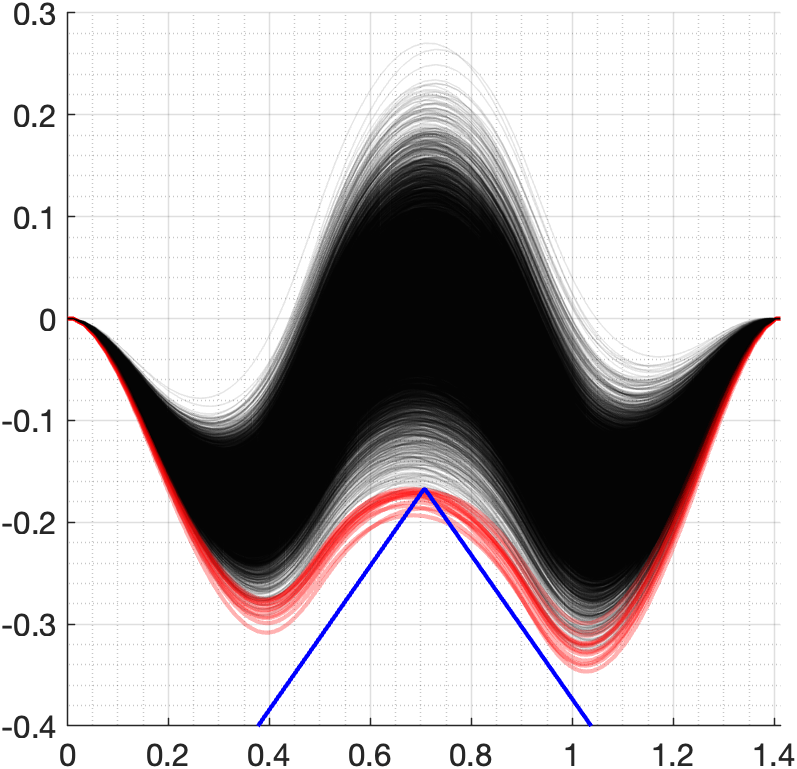} 
         \caption{800 original samples: 14 state constraint violations
                  among 10,000 new samples.}
     \end{subfigure}
     
     \caption{State constraint violations using the optimal control, after
              resampling the normal distribution with 10,000 samples.
              The horizontal axis is the diagonal distance from the origin
              (defining a ``slice'' through the three-dimensional state surfaces).
              The vertical axis represents the interpolated values of the state $u$.
              Shown are the state constraint (blue), the state realizations
              satisfying the constraint (gray), and the state realizations
              violating the constraint (red).}
     \label{fig:stateviolation}
\end{figure}

Fourth, for the case of mean plus CVaR, we evaluate the robustness of the computed optimal
control with respect to the satisfaction of the state inequality constraint, as follows.
We solve the optimal control problem with 100, 200, 400, and 800 samples of the
normal distribution, on a fixed mesh of resolution~64$\times$64.
We then substitute the computed controls into the governing PDE and solve the PDE
for a new set of 10,000 samples of the normal distribution.
Finally, for each of the 10,000 new solutions, we check the pointwise satisfaction
of the state constraint, $u_{\xi} \ge \psi$ over $D$, in the discrete sense and to within
machine tolerance, and report the total number of violations across the 10,000 realizations.
Figure~\ref{fig:stateviolation} demonstrates that the state constraint is violated at
fewer than 1\% of new samples, for all original sample sizes.


\section*{Use of AI tools}
In the final stages of preparing this paper, the authors used the AI models
Claude Opus 5.5 (Anthropic), GPT-6 Astra (OpenAI), and
Gemini 3.1 Pro (Google) to refine proofs and to edit and check the
notation in the text. The AI models were not used to search the literature or to explore or test ideas. All results, arguments and references have been checked by the authors, who
take full responsibility for the content of this paper.

\section*{Funding}
   HA was partially supported by the Office of Naval Research (ONR) under Award NO: N00014-24-1-2147 and 
NSF grant DMS-2408877. HA, DPK and DR were partially supported by the U.S.~Air Force Office of
   Scientific Research, Optimization Program under Award NO: FA9550-25-1-0231. DPK and DR were also supported by the U.S.~Department of Energy's Office of Science Advanced Scientific Computing Research program.
   TMS and HA were supported by SURE-AI: The Norwegian Centre for Sustainable, Risk-Averse, and Ethical AI grant 357482, Research Council of Norway

  Sandia National Laboratories is a multimission laboratory managed and operated by National Technology
       and Engineering Solutions of Sandia, LLC., a wholly owned subsidiary of Honeywell International, Inc.,
       for the U.S. Department of Energy's National Nuclear Security Administration under grant~DE-NA-0003525.
       This paper describes objective technical results and analysis.  Any subjective views or opinions that
       might be expressed in the paper do not necessarily represent the views of the U.S. Department of Energy
       or the United States Government.

\bibliographystyle{siamplain}
\bibliography{refs}  

\end{document}